\documentclass[12pt]{amsart}

\usepackage{amsmath}
\usepackage{amsthm}
\usepackage{amssymb}
\usepackage{latexsym}

\usepackage{eucal}
\usepackage{times}
\usepackage{euler}

\usepackage{pst-node}
\usepackage{pst-plot}
\usepackage{pst-tree}

\usepackage[algoruled,linesnumbered]{algorithm2e}

\author{Piotr \'Sniady}
\address{Institute of Mathematics,
University of Wroclaw,  \mbox{pl.\ Grunwaldzki~2/4,} 50-384
Wroclaw, Poland} \email{Piotr.Sniady@math.uni.wroc.pl}

\title[Generalized Cauchy identities, trees\dots\ Part I]%
{Generalized Cauchy identities, trees \\ and multidimensional
Brownian motions. \\ Part I: bijective proof of generalized Cauchy
identities}

\theoremstyle{plain}
\newtheorem{lemma}{Lemma}
\newtheorem{theorem}[lemma]{Theorem}

\theoremstyle{definition}
\newtheorem*{definition}{Definition}

\theoremstyle{remark}
\newtheorem{remark}[lemma]{Remark}

\newcommand{\A}{{\mathfrak{A}}}

\newcommand{\E}{{\mathbb{E}}}
\newcommand{\C}{{\mathbb{C}}}
\newcommand{\R}{{\mathbb{R}}}

\newcommand{\N}{{\mathbb{N}}}
\newcommand{\gwia}{^{\star}}
\newcommand{\Ha}{{\mathcal{H}}}
\newcommand{\Ka}{{\mathcal{K}}}

\DeclareMathOperator{\tr}{Tr}

\newcommand{\MainPary}{($\alpha$)}
\newcommand{\MainZbiory}{($\beta$)}
\newcommand{\MainCiagi}{($\gamma$)}

\begin{document}

\SetKwFunction{MainBijection}{MainBijection}
\SetKwFunction{SmallBijection}{SmallBijection}
\SetKwFunction{MainBijectionExtended}{MainBijectionExtended}

\begin{abstract}
In this series of articles we study connections between
combinatorics of multidimensional generalizations of the Cauchy
identity and continuous objects such as multidimensional Brownian
motions and Brownian bridges.

In Part I of the series we present a bijective proof of the
multidimensional generalizations of the Cauchy identity. Our
bijection uses oriented planar trees equipped with some linear
orders.

\end{abstract}

\maketitle

\section{Introduction}

Since this paper constitutes the Part I of a series of articles we
allow ourself to start with a longer introduction to the whole
series.

\subsection{Toy example}
\label{subsec:toyexample} The goal of this series of articles is to
discuss multidimensional analogues of the Cauchy identity. However,
before we do this and study our problem in its full generality, we
would like to have a brief look on the simplest case of the (usual)
Cauchy identity. Even in this simplified setting we will be able to
see some important features of the general case.

\subsubsection{Cauchy identity}
\label{subsubsec:Cauchyidentity}

Cauchy identity states that for each nonnegative integer $l$
\begin{equation}
2^{2l}=\sum_{p+q=l} \binom{2p}{p} \binom{2q}{q},
\label{eq:tozsamosc2}
\end{equation}
where the sum runs over nonnegative integers $p,q$. In order to give
a combinatorial meaning to this identity we interpret the left-hand
side of \eqref{eq:tozsamosc2} as the number of sequences
$(x_1,\dots,x_{2l+1})$ such that $x_1,\dots,x_{2l+1}\in\{-1,1\}$ and
$x_1+\dots+x_{2l+1}>0$. For each such a sequence $(x_i)$ we set
$p\geq 0$ to be the biggest integer such that $x_1+\cdots+x_{2p}=0$
and set $q=l-p$; it follows that $(x_i)$ is a concatenation of
sequences $(y_1,\dots,y_{2p})$ and $(z_0,z_1,\dots,z_{2q})$, where
$y_1+\cdots+y_{2p}=0$ and all partial sums of the sequence $(z_i)$
are positive: $z_0+\cdots+z_i>0$ for all $1\leq i\leq 2q$. This can
be illustrated graphically as follows: we treat the sequence $(x_i)$
as a random walk and $2p$ is the time of the last return of the
trajectory to its origin, cf Figure \ref{fig:lamana}. Clearly, for
each value of $p$ there are $\binom{2p}{p}$ ways of choosing the
sequence $(y_i)$ and it is much less obvious (we shall discuss this
problem later on) that for each value of $q$ there are exactly
$\binom{2q}{q}$ ways of choosing the sequence $(z_i)$; in this way
we found a combinatorial interpretation of the right-hand side of
the Cauchy identity \eqref{eq:tozsamosc2}.

\begin{figure}
\psset{unit=0.35cm}
\begin{pspicture}(-0.1,-4.1)(26.3,5.1)
\psaxes[labels=none]{->}(0,0)(-0.1,-4.1)(26.1,5)

\savedata{\mydata}[ {0,-0}, {1,+1}, {2,-0}, {3,+1}, {4,+2}, {5,+3},
{6,+2}, {7,+1}, {8,+2}, {9,+1}, {10,+2}, {11,+1}, {12,-0}, {13,-1},
{14,-2}, {15,-3}, {16,-2}, {17,-1}, {18,0}, {19,+1}, {20,+2},
{21,+3}, {22,+2}, {23,+3}, {24,+2}, {25,+3} ]

\psline[linestyle=dashed,linecolor=gray](18,-4)(18,4)
\psline[linestyle=dashed,linecolor=gray](25,-4)(25,4)

\pcline[linecolor=gray]{<->}(0,-3.5)(18,-3.5) \lput*{:U}{$2p$}
\pcline[linecolor=gray]{<->}(18,-3.5)(25,-3.5) \lput*{:U}{$2q+1$}
\dataplot[showpoints=true]{\mydata}
\end{pspicture}

\caption{A graphical representation of the sequence $(x_i)_{1\leq
i\leq 25}=(1,-1,1,1,1,-1,-1,1,-1,\dots).$ It is also a graph of a
continuous piecewise affine function $X:[0,25]\rightarrow \R$ which
is canonically associated to the sequence $(x_i)$.}
\label{fig:lamana}
\end{figure}

%
%
%
%

\subsubsection{Bijective proof and Pitman transform}
\label{subsubsec:bijectiveproof} In the above discussion we used
without a proof the fact that the number of the sequences
$(z_0,\dots,z_{2q})$ is equal to $\binom{2q}{q}$. The latter number
has a clear combinatorial interpretation as the number of sequences
$(t_1,\dots,t_{2q})$ with $t_1,\dots,t_{2q}\in\{-1,1\}$ and
$t_1+\cdots+t_{2q}=0$, it would be therefore very tempting to proof
the above statement by constructing a bijection between the
sequences $(z_i)$ and the sequences $(t_i)$ and we shall do it in
the following.

Firstly, instead of considering the sequences $(z_0,\dots,z_{2q})$
of length $2q+1$ with all partial sums positive it will be more
convenient to skip the first element and to consider sequences
$(z_1,\dots,z_{2q})$ of length $2q$ such that
$z_1,\dots,z_{2q}\in\{-1,1\}$ with all partial sums nonnegative:
$z_1+\cdots+z_{i}\geq 0$ for all $1\leq i\leq 2q$.

Secondly, it will be convenient to represent the sequences
$(z_1,\dots,z_{2q})$ and $(t_1,\dots,t_{2q})$ as continuous
piecewise affine functions $Z,T:[0,2q]\rightarrow\R$ just as we did
on Figure \ref{fig:lamana}. Formally, function $Z$ is defined as the
unique continuous function such that $Z(0)=0$ and such that for each
integer $1\leq i\leq 2q$ we have $Z'(s)=z_i$ for all $s\in(i-1,i)$.
In this way we can assign a function to any sequence consisting of
only $1$ and $-1$ and we shall make use of this idea later on.

It turns out that an example of a bijection between sequences
$(t_i)$ and $(z_i)$ is provided by the Pitman transform
\cite{Pitman1975} which to a function $T$ associates a function
\begin{equation}
\label{eq:Pitman} Z_s=T_s -2 \inf_{0\leq r\leq s} T_r.
\end{equation}
We shall analyze this bijection in a more general context in Part
III \cite{Sniady2004BijectionAsymptotic} of this series.

\subsubsection{Brownian motion limit and arc-sine law}
\label{subsubsec:Brownianmotionlimit}

What happens to the combinatorial interpretation of the Cauchy
identity \eqref{eq:tozsamosc2} when $l$ tends to infinity? We define
a rescaled function $\tilde{X}_s:[0,1]\rightarrow \R$ given by
$$ \tilde{X}_s = \frac{1}{\sqrt{2l+1}}\ X_{(2l+1) s}, $$
where $X:[0,2l+1]\rightarrow \R$ is the usual function associated to
the sequence $x_1,\dots,x_{2l+1}$ as on Figure \ref{fig:lamana}. The
normalization factors were chosen in such a way that if the sequence
$(x_i)$ is taken randomly (provided $x_1+\cdots+x_{2l+1}>0$) then
the stochastic processes $\tilde{X}_s$ converge  in distribution (as
$l$ tends to infinity) to the Brownian motion $B:[0,1]\rightarrow\R$
conditioned by a requirement that $B_1\geq 0$.

It follows that random variables $\tilde{\Theta}=\sup\ \{t\in[0,1]:
\tilde{X}_t=0\}$ converge in distribution (as $l$ tends to infinity)
to a random variable $\Theta=\sup\ \{ t\in[0,1]: B_t=0\}$, the time
of the last visit of the trajectory of the Brownian motion in the
origin. The discussion from Sections \ref{subsubsec:Cauchyidentity}
and \ref{subsubsec:bijectiveproof} shows that the distribution of
the random variable $\tilde{\Theta}$ is given explicitly by
\begin{equation} \label{eq:exact} P(\tilde{\Theta}<x)=
\sum_{\substack{p+q=l \\ p< x (2l+1) }} \frac{\binom{2p}{p}
\binom{2q}{q}}{2^{2l}}. \end{equation}

The distribution of the random variable $\Theta$ does not change
when we replace the Brownian motion $(B_t)$ conditioned to fulfill
$B_1>0$ by the usual Brownian motion. Thus, Eq.\ \eqref{eq:exact}
after applying the Stirling formula and simple transformations
implies the following well--known result.

\begin{theorem}[Arc-sine law]
If \/ $(B_s)$ is a Brownian motion then the distribution of the
random variable $\Theta=\sup\ \{ t\in[0,1]: B_t=0\}$ is given by
$$ P(\Theta<x)= \frac{1}{2}+\frac{1}{\pi} \sin^{-1}(2x-1) $$
for all $0\leq x\leq 1$.
\end{theorem}

\subsection{How to generalize the Cauchy identity?}
\label{subsec:howtogeneralize}

As we have seen above, the Cauchy identity \eqref{eq:tozsamosc2} has
all properties of a wonderful mathematical result: it is not
obvious, it has interesting applications and it is beautiful. It is
therefore very tempting to look for some more identities which would
share some resemblance to the Cauchy identity or even find some
general identity, equation \eqref{eq:tozsamosc2} would be a special
case of.

Guessing how the left-hand side of \eqref{eq:tozsamosc2} could be
generalized is not difficult and something like $m^{ml}$ is a
reasonable candidate. Unfortunately, it is by no means clear which
sum should replace the right-hand side of \eqref{eq:tozsamosc2}. The
strategy of writing down lots of wild and complicated sums with the
hope of finding the right one by accident is predestined to fail. It
is much more reasonable to find some combinatorial objects which are
counted by the right-hand side of \eqref{eq:tozsamosc2} and then to
find a reasonable generalization of these objects.

 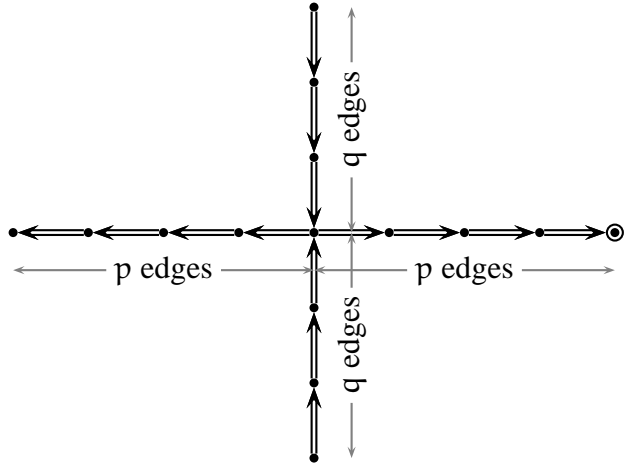
\begin{figure}
 \psset{unit=1cm}
 \begin{pspicture}[](-4,-3.1)(4,3.1)
  \cnode*(4,0){0.6mm}{p1}
  \cnode(4,0){1.2mm}{p1}
  \cnode*(3,0){0.6mm}{p2}
  \cnode*(2,0){0.6mm}{p3}
  \cnode*(1,0){0.6mm}{p4}
  \cnode*(0,0){0.6mm}{cen}

 \ncline[arrowsize=2mm,doubleline=true]{<-}{p1}{p2}
 \ncline[arrowsize=2mm,doubleline=true]{<-}{p2}{p3}
 \ncline[arrowsize=2mm,doubleline=true]{<-}{p3}{p4}
 \ncline[arrowsize=2mm,doubleline=true]{<-}{p4}{cen}

  \cnode*(-4,0){0.6mm}{p1}
  \cnode*(-3,0){0.6mm}{p2}
  \cnode*(-2,0){0.6mm}{p3}
  \cnode*(-1,0){0.6mm}{p4}

 \ncline[arrowsize=2mm,doubleline=true]{<-}{p1}{p2}
 \ncline[arrowsize=2mm,doubleline=true]{<-}{p2}{p3}
 \ncline[arrowsize=2mm,doubleline=true]{<-}{p3}{p4}
 \ncline[arrowsize=2mm,doubleline=true]{<-}{p4}{cen}

  \cnode*(0,-3){0.6mm}{p2}
  \cnode*(0,-2){0.6mm}{p3}
  \cnode*(0,-1){0.6mm}{p4}

 \ncline[arrowsize=2mm,doubleline=true]{->}{p2}{p3}
 \ncline[arrowsize=2mm,doubleline=true]{->}{p3}{p4}
 \ncline[arrowsize=2mm,doubleline=true]{->}{p4}{cen}

  \cnode*(0,3){0.6mm}{p2}
  \cnode*(0,2){0.6mm}{p3}
  \cnode*(0,1){0.6mm}{p4}

 \ncline[arrowsize=2mm,doubleline=true]{->}{p2}{p3}
 \ncline[arrowsize=2mm,doubleline=true]{->}{p3}{p4}
 \ncline[arrowsize=2mm,doubleline=true]{->}{p4}{cen}

 \pcline[linecolor=gray]{<->}(0,-0.5)(4,-0.5) \lput*{:U}{$p$ edges}
 \pcline[linecolor=gray]{<->}(-4,-0.5)(0,-0.5) \lput*{:U}{$p$ edges}
 \pcline[linecolor=gray]{<->}(0.5,0)(0.5,3) \lput*{:U}{$q$ edges}
 \pcline[linecolor=gray]{<->}(0.5,-3)(0.5,0) \lput*{:U}{$q$ edges}

 \end{pspicture}
 \caption{There are $\binom{2p}{p} \binom{2q}{q}$ total orders $<$
 on the vertices of this oriented tree which are compatible with the
 orientation of the edges.}
 \label{fig:trywialnedrzewo}
\end{figure}

For fixed integers $p,q\geq 0$ we consider the tree from Figure
\ref{fig:trywialnedrzewo}. Every edge of this tree is oriented and
it is a good idea to regard these edges as one-way-only roads: if
vertices $x$ and $y$ are connected by an edge and the arrow points
from $y$ to $x$ then the travel from $y$ to $x$ is permitted but the
travel from $x$ to $y$ is not allowed. This orientation defines a
partial order $\prec$ on the set of the vertices: we say that $x
\prec y$ if it is possible to travel from the vertex $y$ to the
vertex $x$ by going through a number of edges (in order to remember
this convention we suggest the Reader to think that $\prec$ is a
simplified arrow $\leftarrow$). Let $<$ be a total order on the set
of the vertices. We say that $<$ is compatible with the orientations
of the edges if for all pairs of vertices $x,y$ such that $x\prec y$
we also have $x<y$. It is very easy to see that for the tree from
Figure \ref{fig:trywialnedrzewo} there are $\binom{2p}{p}
\binom{2q}{q}$ total orders $<$ which are compatible with the
orientations of the edges which coincides with the summand on the
right-hand side of \eqref{eq:tozsamosc2}.

It remains now to find some natural way of generating the trees of
the form depicted on Figure \ref{fig:trywialnedrzewo} with the
property $p+q=l$. We shall do it in the following.

\subsection{Quotient graphs and quotient trees}
\label{subsec:quotientzero}

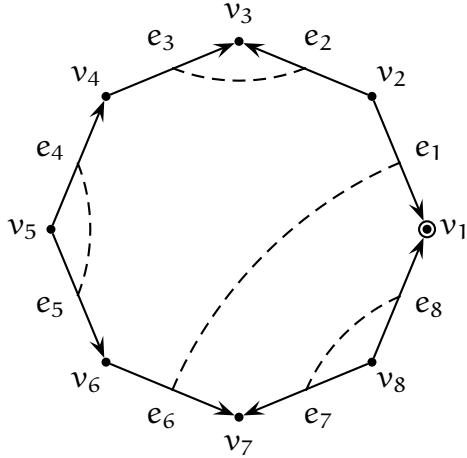
\begin{figure}
\SpecialCoor \psset{unit=0.25cm} \degrees[8]
\begin{pspicture}[](-14,-12)(14,12)
  \cnodeput*(11.5;0){zz}{$v_1$}
 \cnodeput*(11.5;1){zz}{$v_2$}
 \cnodeput*(11.5;2){zz}{$v_3$}
 \cnodeput*(11.5;3){zz}{$v_4$}
 \cnodeput*(11.5;4){zz}{$v_5$}
 \cnodeput*(11.5;5){zz}{$v_6$}
 \cnodeput*(11.5;6){zz}{$v_7$}
 \cnodeput*(11.5;7){zz}{$v_8$}
 \cnodeput*(10.9;0.5){zz}{$e_1$}
 \cnodeput*(10.9;1.5){zz}{$e_2$}
 \cnodeput*(10.9;2.5){zz}{$e_3$}
 \cnodeput*(10.9;3.5){zz}{$e_4$}
 \cnodeput*(10.9;4.5){zz}{$e_5$}
 \cnodeput*(10.9;5.5){zz}{$e_6$}
 \cnodeput*(10.9;6.5){zz}{$e_7$}
 \cnodeput*(10.9;7.5){zz}{$e_8$}
 \cnode*(10;0){0.6mm}{p1}
  \cnode(10;0){1.2mm}{p1}
 \cnode*(10;1){0.6mm}{p2}
 \cnode*(10;2){0.6mm}{p3}
 \cnode*(10;3){0.6mm}{p4}
 \cnode*(10;4){0.6mm}{p5}
 \cnode*(10;5){0.6mm}{p6}
 \cnode*(10;6){0.6mm}{p7}
 \cnode*(10;7){0.6mm}{p8}
 \ncline[arrowsize=2mm]{<-}{p1}{p2}
 \ncline[arrowsize=2mm]{->}{p2}{p3}
 \ncline[arrowsize=2mm]{<-}{p3}{p4}
 \ncline[arrowsize=2mm]{<-}{p4}{p5}
 \ncline[arrowsize=2mm]{->}{p5}{p6}
 \ncline[arrowsize=2mm]{->}{p6}{p7}
 \ncline[arrowsize=2mm]{<-}{p7}{p8}
 \ncline[arrowsize=2mm]{->}{p8}{p1}
 \pnode(9.3;0.5){c1}
 \pnode(9.3;1.5){c2}
 \pnode(9.3;2.5){c3}
 \pnode(9.3;3.5){c4}
 \pnode(9.3;4.5){c5}
 \pnode(9.3;5.5){c6}
 \pnode(9.3;6.5){c7}
 \pnode(9.3;7.5){c8}
 \ncarc[linestyle=dashed,arcangle=-0.5]{c1}{c6}
 \ncarc[linestyle=dashed,arcangle=0.5]{c2}{c3}
 \ncarc[linestyle=dashed,arcangle=0.5]{c4}{c5}
 \ncarc[linestyle=dashed,arcangle=0.5]{c7}{c8}

\end{pspicture}
\caption{A graph $G_{\epsilon}$ corresponding to the sequence
$\epsilon=(+1,-1,+1,+1,-1,-1,+1,-1)$. The dashed lines represent the
pairing $\sigma=\big\{ \{1,6\},\{2,3\},\{4,5\},\{7,8\} \} \big\}$.}
\label{fig:graphG}
\end{figure}

We recall now the construction of Dykema and Haagerup
\cite{DykemaHaagerup2001}. For integer $k\geq 1$ let $G$ be an
oriented $k$--gon graph with consecutive vertices $v_1,\dots,v_k$
and edges $e_1,\dots,e_k$ (edge $e_i$ connects vertices $v_i$ and
$v_{i+1}$). The vertex $v_1$ is distinguished, see Figure
\ref{fig:graphG}. We encode the information about the orientations
of the edges in a sequence $\epsilon(1),\dots,\epsilon(k)$ where
$\epsilon(i)=+1$ if the arrow points from $v_{i+1}$ to $v_i$ and
$\epsilon(i)=-1$ if the arrow points from $v_{i}$ to $v_{i+1}$. The
graph $G$ is uniquely determined by the sequence $\epsilon$ and
sometimes we will explicitly state this dependence by using the
notation $G_{\epsilon}$.

\begin{figure}
\SpecialCoor \psset{unit=0.6cm} \degrees[360]
 \begin{pspicture}[](-5,-5)(14,5)
 \cnode*(0;0){0.6mm}{A}
 \cnode*(5;0){0.6mm}{B}
 \cnode(5;0){1.2mm}{B}
 \cnode*(10;0){0.6mm}{C}
 \cnode*(5;120){0.6mm}{D}
 \cnode*(5;240){0.6mm}{E}
 \ncline[arrowsize=2mm,doubleline=true]{->}{A}{B}
 \ncline[arrowsize=2mm,doubleline=true]{<-}{B}{C}
 \ncline[arrowsize=2mm,doubleline=true]{->}{A}{D}
 \ncline[arrowsize=2mm,doubleline=true]{<-}{A}{E}
 \rput[b](2.5,0.2){$e_1$}
 \rput[t](2.5,-0.2){$e_6$}
 \rput[b](7.5,0.2){$e_8$}
 \rput[t](7.5,-0.2){$e_7$}
 \rput[bl](2.5;115){$e_2$}
 \rput[tr](2.5;125){$e_3$}
 \rput[br](2.5;235){$e_4$}
 \rput[tl](2.5;245){$e_5$}
 \rput[t](5,-0.4){$R=v_1=v_7$}
 \rput[t](10,-0.4){$v_8$}
 \rput[r](-0.3,0){$v_2=v_4=v_6$}
 \rput[br](5.3;120){$v_3$}
 \rput[tr](5.3;240){$v_5$}
 \end{pspicture}
 \caption{The quotient graph $T_{\sigma}$ corresponding to the graph from Figure
 \ref{fig:graphG}. The root $R$ of the tree $T_{\sigma}$ is encircled.}
 \label{fig:quotient}
\end{figure}
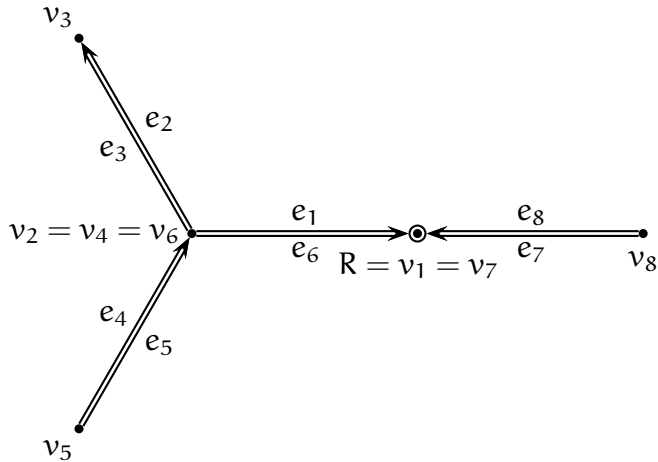

Let $\sigma=\big\{ \{i_1,j_1\},\dots,\{i_{k/2},j_{k/2}\} \big\}$ be
a pairing of the set $\{1,\dots,k\}$, i.e.\ pairs $\{i_m,j_m\}$ are
disjoint and their union is equal to $\{1,\dots,k\}$. We say that
$\sigma$ is compatible with $\epsilon$ if
\begin{equation} \epsilon(i)+\epsilon(j)=0 \qquad \text{ for every }\{i,j\}\in\sigma.
\label{eq:orientacja} \end{equation} It is a good idea to think that
$\sigma$ is a pairing between the edges of $G$, see Figure
\ref{fig:graphG}. For each $\{i,j\}\in\sigma$ we identify (or, in
other words, we glue together) the edges $e_i$ and $e_j$ in such a
way that the vertex $v_i$ is identified with $v_{j+1}$ and vertex
$v_{i+1}$ is identified with $v_{j}$ and we denote by $T_{\sigma}$
the resulting quotient graph. Since each edge of $T_{\sigma}$
origins from a pair of edges of $G$, we draw all edges of
$T_{\sigma}$ as double lines. The condition \eqref{eq:orientacja}
implies that each edge of $T_{\sigma}$ carries a natural
orientation, inherited from each of the two edges of $G$ it comes
from, see Figure \ref{fig:quotient}.

From the following on, we consider only the case when the quotient
graph $T_{\sigma}$ is a tree. One can show \cite{DykemaHaagerup2001}
that the latter holds if and only if the pairing $\sigma$ is
non--crossing \cite{Kreweras}; in other words it is not possible
that for some $p<q<r<s$ we have $\{p,r\},\{q,s\}\in\sigma$. The name
of the non--crossing pairings comes from their property that on
their graphical depictions (such as Figure \ref{fig:graphG}) the
lines do not cross. Let the root $R$ of the tree $T_{\sigma}$ be the
vertex corresponding to the distinguished vertex $v_1$ of the graph
$G$.

\subsection{How to generalize the Cauchy identity? (continued)}
\label{subsec:howtogeneralize2} 
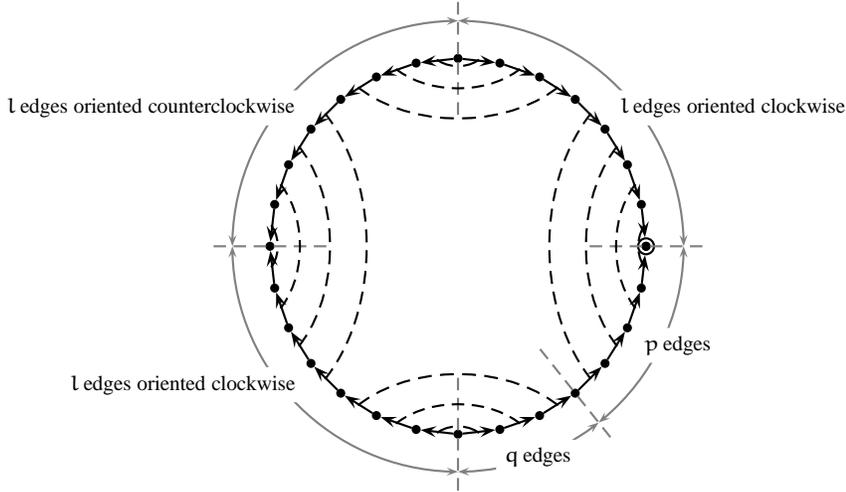
\begin{figure}
\SpecialCoor \psset{unit=0.25cm} \degrees[28]
\begin{pspicture}[](-14,-12)(14,12)
 \cnode*(10;0){0.6mm}{p1}
 \cnode(10;0){1.2mm}{p1}

 \cnode*(10;1){0.6mm}{p2}
 \cnode*(10;2){0.6mm}{p3}
 \cnode*(10;3){0.6mm}{p4}
 \cnode*(10;4){0.6mm}{p5}
 \cnode*(10;5){0.6mm}{p6}
 \cnode*(10;6){0.6mm}{p7}
 \cnode*(10;7){0.6mm}{p8}
 \cnode*(10;8){0.6mm}{p9}
 \cnode*(10;9){0.6mm}{p10}

 \cnode*(10;10){0.6mm}{p11}
 \cnode*(10;11){0.6mm}{p12}
 \cnode*(10;12){0.6mm}{p13}
 \cnode*(10;13){0.6mm}{p14}
 \cnode*(10;14){0.6mm}{p15}
 \cnode*(10;15){0.6mm}{p16}
 \cnode*(10;16){0.6mm}{p17}
 \cnode*(10;17){0.6mm}{p18}
 \cnode*(10;18){0.6mm}{p19}
 \cnode*(10;19){0.6mm}{p20}

 \cnode*(10;20){0.6mm}{p21}
 \cnode*(10;21){0.6mm}{p22}
 \cnode*(10;22){0.6mm}{p23}
 \cnode*(10;23){0.6mm}{p24}
 \cnode*(10;24){0.6mm}{p25}
 \cnode*(10;25){0.6mm}{p26}
 \cnode*(10;26){0.6mm}{p27}
 \cnode*(10;27){0.6mm}{p28}
 \cnode*(10;28){0.6mm}{p29}

 \ncline[arrowsize=1.5mm]{<-}{p1}{p2}
 \ncline[arrowsize=1.5mm]{<-}{p2}{p3}
 \ncline[arrowsize=1.5mm]{<-}{p3}{p4}
 \ncline[arrowsize=1.5mm]{<-}{p4}{p5}
 \ncline[arrowsize=1.5mm]{<-}{p5}{p6}
 \ncline[arrowsize=1.5mm]{<-}{p6}{p7}
 \ncline[arrowsize=1.5mm]{<-}{p7}{p8}

 \ncline[arrowsize=1.5mm]{<-}{p15}{p16}
 \ncline[arrowsize=1.5mm]{<-}{p16}{p17}
 \ncline[arrowsize=1.5mm]{<-}{p17}{p18}
 \ncline[arrowsize=1.5mm]{<-}{p18}{p19}
 \ncline[arrowsize=1.5mm]{<-}{p19}{p20}
 \ncline[arrowsize=1.5mm]{<-}{p20}{p21}
 \ncline[arrowsize=1.5mm]{<-}{p21}{p22}

 \ncline[arrowsize=1.5mm]{->}{p8}{p9}
 \ncline[arrowsize=1.5mm]{->}{p9}{p10}
 \ncline[arrowsize=1.5mm]{->}{p10}{p11}
 \ncline[arrowsize=1.5mm]{->}{p11}{p12}
 \ncline[arrowsize=1.5mm]{->}{p12}{p13}
 \ncline[arrowsize=1.5mm]{->}{p13}{p14}
 \ncline[arrowsize=1.5mm]{->}{p14}{p15}

 \ncline[arrowsize=1.5mm]{->}{p22}{p23}
 \ncline[arrowsize=1.5mm]{->}{p23}{p24}
 \ncline[arrowsize=1.5mm]{->}{p24}{p25}
 \ncline[arrowsize=1.5mm]{->}{p25}{p26}
 \ncline[arrowsize=1.5mm]{->}{p26}{p27}
 \ncline[arrowsize=1.5mm]{->}{p27}{p28}
 \ncline[arrowsize=1.5mm]{->}{p28}{p1}

 \psline[linecolor=gray,linestyle=dashed](7;0)(13;0)
 \psline[linecolor=gray,linestyle=dashed](7;7)(13;7)
 \psline[linecolor=gray,linestyle=dashed](7;14)(13;14)
 \psline[linecolor=gray,linestyle=dashed](7;21)(13;21)
 \psline[linecolor=gray,linestyle=dashed](7;24)(13;24)

 \psarc[linecolor=gray]{<->}(0,0){12}{0}{7}
 \psarc[linecolor=gray]{<->}(0,0){12}{7}{14}
 \psarc[linecolor=gray]{<->}(0,0){12}{14}{21}
 \psarc[linecolor=gray]{<->}(0,0){12}{21}{24}
 \psarc[linecolor=gray]{<->}(0,0){12}{24}{28}

 \rput*[bl]{0}(11;3){{\tiny $l$ edges oriented clockwise}}
 \rput*[br]{0}(11;11){{\tiny $l$ edges oriented counterclockwise}}
 \rput*[tr]{0}(11;17){{\tiny $l$ edges oriented clockwise}}
 \rput*[tl]{0}(11;26){{\tiny $p$ edges}}
 \rput*[tl]{0}(11;22){{\tiny $q$ edges}}

 \pcarc[linestyle=dashed,arcangle=-3](10;0.5)(10;27.5)
 \pcarc[linestyle=dashed,arcangle=-3](10;1.5)(10;26.5)
 \pcarc[linestyle=dashed,arcangle=-3](10;2.5)(10;25.5)
 \pcarc[linestyle=dashed,arcangle=-3](10;3.5)(10;24.5)
 \pcarc[linestyle=dashed,arcangle=-3](10;14.5)(10;13.5)
 \pcarc[linestyle=dashed,arcangle=-3](10;15.5)(10;12.5)
 \pcarc[linestyle=dashed,arcangle=-3](10;16.5)(10;11.5)
 \pcarc[linestyle=dashed,arcangle=-3](10;17.5)(10;10.5)
 \pcarc[linestyle=dashed,arcangle=-3](10;7.5)(10;6.5)
 \pcarc[linestyle=dashed,arcangle=-3](10;8.5)(10;5.5)
 \pcarc[linestyle=dashed,arcangle=-3](10;9.5)(10;4.5)
 \pcarc[linestyle=dashed,arcangle=-3](10;21.5)(10;20.5)
 \pcarc[linestyle=dashed,arcangle=-3](10;22.5)(10;19.5)
 \pcarc[linestyle=dashed,arcangle=-3](10;23.5)(10;18.5)

\end{pspicture}
\caption[ssss]{A graph $T$ corresponding to sequence
$\epsilon=(\underbrace{+1}_{l \text{ times}}, \underbrace{-1}_{l
\text{ times}}, \underbrace{+1}_{l \text{ times}},
\underbrace{-1}_{l \text{ times}})$. The dashed lines denote a
pairing $\sigma$ for which the quotient graph $T_{\sigma}$ is
depicted on Figure \ref{fig:trywialnedrzewo}.}

\label{fig:graphkwadratowy}
\end{figure}

Let us come back to the discussion
from Section \ref{subsec:howtogeneralize}. We consider the polygon
$G_{\epsilon}$ corresponding to
$$\epsilon=(\underbrace{+1}_{l \text{ times}},
\underbrace{-1}_{l \text{ times}}, \underbrace{+1}_{l \text{
times}}, \underbrace{-1}_{l \text{ times}}).$$ All possible
non-crossing pairings $\sigma$ which are compatible with $\epsilon$
are depicted on Figure \ref{fig:graphkwadratowy} and it easy to see
that the corresponding quotient tree $T_{\sigma}$ has exactly the
form depicted on Figure \ref{fig:trywialnedrzewo}.

In this way we managed to find relatively natural combinatorial
objects, the number of which is given by the right-hand side of the
Cauchy identity \eqref{eq:tozsamosc2}. After some guesswork we end
up with the following conjecture (please note that the usual Cauchy
identity \eqref{eq:tozsamosc2} corresponds to $m=2$).

\begin{theorem}[Generalized Cauchy identity]
\label{theo:generalized} For integers $l,m\geq 1$ there are exactly
$m^{ml}$ pairs $(\sigma,<)$, where $\sigma$ is a non-crossing
pairing compatible with \begin{equation} \label{eq:epsilonnk}
\epsilon=(\underbrace{\underbrace{+1}_{l \text{
times}},\underbrace{-1}_{l \text{ times}}, \underbrace{+1}_{l \text{
times}},\underbrace{-1}_{l \text{ times}},\dots}_{2m \text{ blocks,
i.e.\ total of\/  } 2ml \text{ elements}})
\end{equation}
and $<$ is a total order on the vertices of $T_{\sigma}$ which is
compatible with the orientations of the edges.
\end{theorem}

Above we provided only vague heuristical arguments why the above
conjecture could be true. Surprisingly, as we shall see in the
following, Theorem \ref{theo:generalized} is indeed true.

The formulation of Theorem \ref{theo:generalized} is combinatorial
and therefore appears to be far from its motivation, the usual
Cauchy identity \eqref{eq:tozsamosc2}, which is formulated
algebraically, nevertheless for each fixed value of $m$ one can
enumerate all `classes' of pairings compatible with
\eqref{eq:epsilonnk} and for each class count the number of
compatible orders $<$. To give to the Reader a flavor of the
algebraic implications of Theorem \ref{theo:generalized}, we present
the case of $m=3$ \cite{DykemaYan}
\begin{multline}
\label{eq:tozsamosc3}
3^{3l} =\sum_{p+q=l} \binom{3p}{p,p,p} \binom{3q}{q,q,q} + \\
 + 3 \sum_{\substack{p+q+r=l-1\\ r'+q'=r+q+1 \\ p''+r''=p+r+1}}
\binom{2p+p''}{p,p,p''} \binom{2q+q'}{q,q,q'}
\binom{r+r'+r''}{r,r',r''}.
\end{multline}
The complication of the formula grows very quickly and already for
$m=4$ the appropriate expression has a length of a half page of a
typed text \cite{Sniady2002}.

%

\subsection{Historical overview: operator algebras,
     free probability and triangular operator $T$}
The history presented in Sections
\ref{subsec:howtogeneralize}--\ref{subsec:howtogeneralize2} of
finding the generalization of the Cauchy identity is too nice to be
true and indeed it is not the way how Theorem \ref{theo:generalized}
was postulated. As we shall see in the following, the path towards
this result led not through combinatorics but through the theory of
operator algebras. Since this section is very loosely connected with
the rest of this article, Readers not interested in theory operator
algebras may skip it without much harm.

\subsubsection{Invariant subspace conjecture}
One of the fundamental problems of the theory of operator algebras
is the invariant subspace conjecture which asks if for every bounded
operator $x$ acting on an infinite-dimensional Hilbert space $\Ha$
there exists a closed subspace $\Ka\subset\Ha$ such that $\Ka$ is
nontrivial in the sense that $\Ka\neq \{0\}$, $\Ka\neq \Ha$ and
which is an invariant subspace of $x$. Since for many decades nobody
was able to prove the invariant subspace conjecture in its full
generality, Dykema and Haagerup took the opposite strategy and tried
to construct explicitly a counterexample by the means of the
Voiculescu's free probability theory.

The free probability \cite{VoiculescuDykemaNica,HiaiPetz2000} is a
non-commutative probability theory with the classical notion of
independence replaced by the notion of freeness. Natural examples
which fit nicely into the framework of the free probability include
large random matrices, free products of von Neumann algebras and
asymptotics of large Young diagrams. Families of operators which
arise in the free probability are, informally speaking, very
non-commutative and for this reason they are perfect candidates for
counterexamples to the conjectures in the theory of operator
algebras \cite{VoiculescuPart3}.

The first candidate for a counterexample to the invariant subspace
conjecture considered by Dykema and Haagerup was the circular
operator, which unfortunately turned out to have a large family of
invariant subspaces \cite{DykemaHaagerup2000,SniadySpeicher2001}.
Later on Haagerup \cite{Haagerup2001} proved a version of a spectral
theorem for certain non--normal operators and thus he constructed
invariant subspaces for many classes of operators. This result gave
very strong restrictions on the form of a possible counterexample,
namely the Brown spectral measure \cite{Brown} of such an operator
should be concentrated in only one point. It was a hint to look for
counterexamples among, so-called, quasinilpotent operators. In this
way Dykema and Haagerup \cite{DykemaHaagerup2001} initiated a study
of the triangular operator $T$, which appeared at that time to be a
perfect candidate because it is quasinilpotent and it admits very
nice random matrix models.

\subsubsection{Triangular operator $T$}
The triangular operator $T$ \cite{DykemaHaagerup2001} can be
abstractly described as an element of a von Neumann algebra $\A$
equipped with a finite normal faithful tracial state
$\phi:\A\rightarrow\C$ with the non-commutative moments
$\phi(T^{\epsilon(1)} \cdots T^{\epsilon(n)})$ given by
\begin{equation}
\phi(T^{\epsilon(1)} \cdots
T^{\epsilon(n)})=\lim_{N\rightarrow\infty} \frac{1}{N} \E \tr
(T_N^{\epsilon(1)} \cdots T_N^{\epsilon(n)})
\label{eq:oryginalnadefinicja}
\end{equation}
for any $n\in\N$ and $\epsilon(1),\dots,\epsilon(n)\in\{-1,+1\}$,
where we use the notation $T^{+1}:=T$ and $T^{-1}:=T\gwia$; and
where
$$T_N = \left[ \begin{array}{ccccc} t_{1,1} &
t_{1,2} &   \cdots & t_{1,n-1} & t_{1,n} \\ 0 & t_{2,2} & \cdots &
t_{2,n-1} & t_{2,n} \\ \vdots&                 & \ddots & \vdots   &
\vdots        \\ & & & t_{n-1,n-1} & t_{n-1,n} \\ 0&       & \cdots
&             0               & t_{n,n} \end{array} \right]$$ is an
upper-triangular random matrix, the entries $(t_{i,j})_{1\leq i\leq
j\leq N}$ of which are independent centered Gaussian random
variables with variance $\frac{1}{N}$.

Definition \eqref{eq:oryginalnadefinicja} is not very convenient and
one can show \cite{DykemaHaagerup2001} that it is equivalent to the
following one:  $(n/2 +1)!\ \phi(T^{\epsilon(1)} \cdots
T^{\epsilon(n)})$ is equal to the number of pairs $(\sigma,<)$ such
that $\sigma$ is a pairing compatible with $\epsilon$ and $<$ is a
total order on the vertices of $T_{\sigma}$ which is compatible with
the orientation of the edges. The Reader may easily see that the
latter definition of $T$ is very closely related to the results
presented in this paper; in particular Theorem
\ref{theo:generalized} can be now equivalently stated as follows (in
fact it is the form in which Dykema and Haagerup stated originally
their conjecture \cite{DykemaHaagerup2001}):
\begin{theorem}
\label{theo:momentybyly} If \/ $l,m\geq 1$ are integers then
$$ \phi \left[\big(T^l (T\gwia)^l\big)^m \right] = \frac{m^{ml}}{(ml+1)!}. $$
\end{theorem}

Yet another approach to $T$ is connected with the combinatorial
approach to operator-valued free probability \cite{Speicher1998},
namely $T$ can be described as a certain generalized circular
element. Speaking very briefly, the non-commutative moments of $T$
can be described as certain iterated integrals \cite{Sniady2002}.
This approach turned out to be very fruitful: in this way in our
previous work \cite{Sniady2002} we found the first proof of Theorem
\ref{theo:generalized} and Theorem \ref{theo:momentybyly}; a
different proof was later presented in \cite{AagaardHaagerup2004}.
Some other combinatorial results concerning the non-commutative
moments of $T$ were obtained in \cite{DykemaYan}.

Theorem \ref{theo:generalized} and Theorem \ref{theo:momentybyly}
were conjectured by Dykema and Haagerup \cite{DykemaHaagerup2001} in
the hope that they might be useful in the study of spectral
properties of $T$. Literally speaking, this hope turned out to be
wrong since the later construction of the hyperinvariant subspaces
of $T$ by Dykema and Haagerup \cite{DykemaHaagerup2003,
Haagerup2002Beijing} did not make use of Theorem
\ref{theo:generalized} and Theorem \ref{theo:momentybyly}, however
it made use of one of the auxiliary results used in our proof
\cite{Sniady2002} of these theorems. In this way, indirectly,
Theorem \ref{theo:generalized} and Theorem \ref{theo:momentybyly}
turned out to be indeed helpful for their original purpose.

As we already mentioned, Dykema and Haagerup
\cite{DykemaHaagerup2003, Haagerup2002Beijing} constructed a family
of hyperinvariant subspaces of $T$ and in this way the original
motivation for studying the operator $T$ (as a possible
counterexample for the invariant subspace conjecture) ended up as a
failure. There are still some investigations of the triangular
operator $T$ as a possible counterexample for some other
conjectures, for example \cite{Aagaard2003}, however most of the
specialists do not expect any surprises in the theory of operator
algebras coming from this direction. In this article we would like
to convince the Reader that the applications of the triangular
operator $T$ in combinatorics and the classical probability theory
constitute a sufficient compensation for the lost hopes concerning
its applications in the theory of operator algebras.

\subsection{Overview of this series of articles}
\nopagebreak
\subsubsection{Part I: Bijective proof of generalized Cauchy
identities }

\label{subsubsec:ocochodziwpierwszejczesci} In this article we shall
prove the following result. Let $l_1\leq \cdots \leq l_m$ be a
weakly increasing sequence of positive integers; we denote
$L=l_1+\cdots+l_m+1$.
For $1\leq i\leq m$ we set
\begin{equation} \epsilon_i= \big(
\underbrace{+1}_{l_i \text{ times}},\underbrace{-1}_{l_{i-1} \text{ times}},
\dots,\underbrace{(-1)^{i-1}}_{l_1 \text{
times}}, \underbrace{(-1)^{i}}_{l_1 \text{ times}}, \dots,
\underbrace{+1}_{l_{i-1} \text{ times}},\underbrace{-1}_{l_{i} \text{ times}}
\big),
\label{eq:epsilon}
\end{equation}
where $\underbrace{a}_{l \text{ times}}$ denotes
$\underbrace{a,\dots,a}_{l \text{ times}}$. The Reader may restrict
his/her attention to the most interesting case when
$l_1=l_2=\cdots=l$ are all equal and $\epsilon_i$ takes a simpler
form
$$ \epsilon_i= \big(  \underbrace{\underbrace{+1}_{l \text{ times}},
\underbrace{-1}_{l \text{ times}}}_{i \text{ times, i.e.\ a
total of $2li$ elements}} \big).$$ In this case $\epsilon_m$ 
coincides with \eqref{eq:epsilonnk} and it is easy to show from the Raney
lemma \cite{Raney1960} that
the set {\MainZbiory} below has $m^{ml}$ elements \cite{Sniady2002}
hence Theorem \ref{theo:generalized} will follow from the following
stronger result.

\begin{theorem}[The main result]
\label{theo:mainresult} Let $\epsilon_{m}$ 
be as above. The algorithm \MainBijection described in this article provides a
bijection between
\begin{enumerate}
\item[\MainPary] the set of pairs $(\sigma,<)$, where $\sigma$ is a
pairing compatible with $\epsilon_{m}$ and $<$ is a total order on
the vertices of\/ $T_{\sigma}$ which is compatible with the
orientations of the edges; \label{bij5}

\item[\MainZbiory] the set of tuples $(B_1,\dots,B_m)$, where
$B_1,\dots,B_m$ are disjoint sets such that $B_1 \cup \cdots \cup
B_m=\{1,\dots,L\}$ and
$$ |B_1|+\cdots+|B_n|\leq l_1+\cdots+l_n$$
holds true for each $1\leq n \leq m-1$;
\end{enumerate}

Alternatively, set {\MainZbiory} can be described as
\begin{enumerate}
\item[\MainCiagi] the set of sequences $(a_1,\dots,a_{L})$ such that
$a_1,\dots,a_{L}\in\{1,\dots,m\}$ and for each $1\leq n \leq m-1$ at
most $l_1+\cdots+l_n$ elements of the sequence $(a_i)$ belong to the
set $\{1,\dots,n\}$;
\end{enumerate}
where the bijection between sets {\MainZbiory} and {\MainCiagi} is
given by $B_j=\{k: a_k=j\}$. 
\end{theorem}

\begin{remark}
The sequence $(\tilde{a}_1,\dots,\tilde{a}_{L})$ can be regarded as a
generalized parking function, where $\tilde{a}_r=m+1 - a_r$. Indeed, let
$(b_1,\dots,b_L)$ be its non-decreasing rearrangement; then the original
sequence $(a_1,\dots,a_L)$ contributes to {\MainCiagi} iff $b_1,\dots,b_L$ are
positive integers such that
$b_{1+l_m}\leq 1$,
$b_{1+l_m+l_{m-1}}\leq 2$,\dots, $b_{1+l_m+\cdots+l_1}\leq m$ which is a
slighlty modified definition of a parking function.
\end{remark}

The bijection provided by the above theorem plays the central role
in this series of articles.


\subsubsection{Part II: Combinatorial differential calculus
\cite{Sniady2004BijectionDifferential}}

The bijections considered in Part I of this series (Section
\ref{sec:almostthemain} and Section \ref{sec:themainresult} of this
article) are far from being trivial and the Reader might wonder how
did the author guess their correct form and what is the conceptual
idea behind them. To answer these questions we would like to come
back to our previous work \cite{Sniady2002} where we provided the
first proof of Theorem \ref{theo:generalized}. The main idea was to
associate a polynomial of a single variable to every pair
$(\sigma,<)$ and by additivity to every graph $G_{\epsilon}$. The
polynomials associated to $\epsilon$ as in \eqref{eq:epsilonnk} with
different values of $m$ turned out to be related by a simple
differential equation and for this reason can be regarded as
generalizations of Abel polynomials.

In Part II of this series \cite{Sniady2004BijectionDifferential}
(joint work with Artur Je\.z) we present an analogue of the
differential calculus in which the role of polynomials is played by
certain ordered sets and trees. Our combinatorial calculus has all
nice features of the usual calculus and has an advantage that the
elements of the considered ordered sets might carry some additional
information. In this way our analytic proof from \cite{Sniady2002}
can be directly reformulated in our new language of the
combinatorial calculus; furthermore the additional information
carried by the vertices determines uniquely the bijections presented
in Part I of this series.


\subsubsection{Part III: Multidimensional arc-sine laws
\cite{Sniady2004BijectionAsymptotic}}
In Section \ref{subsubsec:Brownianmotionlimit} we presented how a
bijective proof the usual Cauchy identity can be used to extract
some information about the behavior of the Brownian motion and in
particular to show the arc-sine law. It is therefore natural to ask
if the bijective proof of the generalized Cauchy identities
presented in Part I could provide some information about
multidimensional Brownian motions.

In order to answer these questions we study in Part III of this
series the asymptotic behavior of the trees and bijections presented
in Part I. Asymptotically, as their size tend to infinity, these
trees converge towards continuous objects such as multidimensional
Brownian motions and Brownian bridges. Our bijection behaves nicely
in this asymptotic setting and becomes a map between certain classes
of functions valued in $\R^{m-1}$, which is closely related to the
Pitman transform and Littelmann paths. In this way we are able to
describe certain interesting properties of multidimensional Brownian
motions and in particular we prove a multidimensional analogue of
the arc-sine law.

%




%
%


\section{The main bijection}

\subsection{Structure of a planar tree. Order $\lhd$}
\label{sec:preorder} For a non--crossing pairing $\sigma$ we can
describe the process of creating the quotient graph as follows: we
think that the edges of the graph $G$ are sticks of equal lengths
with flexible connections at the vertices. Graph $G$ is lying on a
flat surface in such a way that the edges do not cross. For each
pair $\{i,j\}\in\sigma$ we glue together edges $e_i$ and $e_j$ by
bending the joints in such a way that the sticks should not cross.
In this way $T_{\sigma}$ has a structure of a planar tree, i.e.\ for
each vertex we can order the adjacent edges up to a cyclic shift
(just like points on a circle). We shall provide an alternative
description of this planar structure in the following.

Let us visit the vertices of $G$ in the usual cyclic order
$v_1,v_2,\dots,v_k,v_1$ by going along the edges $e_1,\dots,e_k$; by
passing to the quotient graph $T_{\sigma}$ we obtain a journey on
the graph $T_{\sigma}$ which starts and ends in the root $R$. The
structure of the planar tree defined above can be described as
follows: if we travel on the graphical representation of
$T_{\sigma}$ by touching the edges by our left hand, we obtain the
same journey. For each vertex of $T_{\sigma}$ we mark the time we
visit it for the first time; comparison of these times gives us a
total order $\lhd$, called preorder \cite{StanleyEnumerativeVol2},
on the vertices of $T_{\sigma}$. For example, in the case of the
tree from Figure \ref{fig:quotient} we have $v_1\lhd v_2\lhd v_3\lhd
v_5\lhd v_8$.

\subsection{Pairing between leafs and bays}
\label{subsec:leafsandbays}

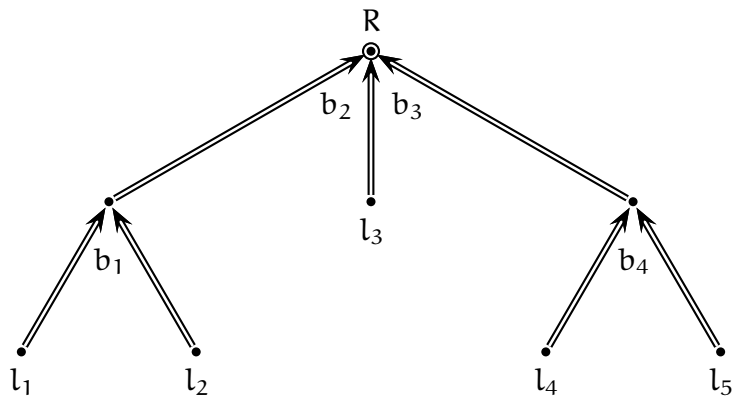
\begin{figure}
\renewcommand{\psedge}[2]{\ncline[doubleline=true]{<-}{#1}{#2}}
\psset{unit=1cm}
\begin{pspicture}[](0,-4.5)(0,0.5)
\pstree[nodesep=0.3mm,treemode=D,treefit=loose,treesep=2]{
\Tc{1.2mm}~[tnpos=a,tnsep=0mm]{$R$} \Tc*{0.6mm} }{
\pstree{\Tc*{0.6mm}~[tnpos=b,tnsep=5mm]{$b_1$}}{
\Tc*{0.6mm}~[tnpos=b,tnsep=1mm]{$l_1$}
\Tc*{0.6mm}~[tnpos=b,tnsep=1mm]{$l_2$} }
\Tc*{0.6mm}~[tnpos=b,tnsep=1mm]{$l_3$}  \tlput[tpos=0.35]{$b_2$}
\trput[tpos=0.35]{$b_3$}
\pstree{\Tc*{0.6mm}~[tnpos=b,tnsep=5mm]{$b_4$}}{
\Tc*{0.6mm}~[tnpos=b,tnsep=1mm]{$l_4$}
\Tc*{0.6mm}~[tnpos=b,tnsep=1mm]{$l_5$}}}
\end{pspicture}
\caption{Example of a tree such that the arrows on all the edges
point towards the root. Leafs $l_1,l_2,\dots$ and bays
$b_1,b_2,\dots$ are marked.}
 \label{fig:leafsandbays}
\end{figure}

Suppose that $U$ is an oriented planar tree with the property that
the arrows on all the edges are pointing towards the root $R$; in
other words $R\preceq x$ holds true for every vertex $x$. We shall
also assume that the tree $U$ consists of at least two vertices.

We call a pair of edges $\{e,f\}$ a bay if edges $e,f$ share a
common vertex $v$ and are adjacent edges (adjacent with respect to
the structure of the planar tree) and arrows on $e$ and $f$ point
towards the common vertex $v$. It is convenient to represent a bay
as the corner between edges $e$ and $f$, cf Figure
\ref{fig:leafsandbays}.

A vertex is called a leaf if it is connected to exactly one edge and
it is different from the root $R$, cf Figure \ref{fig:leafsandbays}.

Let us travel on the tree $U$ (we begin and end at the root $R$) in
such a way that we always touch the edges of the tree by our left
hand. We say that a passage along an edge is negative if the arrow
on the edge coincides with the direction of travel; otherwise we
call it a positive passage (the origin of this convention is the
following: if $U=T_{\sigma}$ is a quotient tree coming from a
polygonal graph $G_{\epsilon}$, where
$\epsilon=(\epsilon(1),\dots,\epsilon(k))$ then the sign of the
$n$-th step coincides with the sign of $\epsilon(n)$). It is easy to
see that a bay corresponds to a pair of consecutive passages: a
negative and a positive one; similarly entering and leaving a leaf
corresponds to a pair of consecutive passages: a positive and a
negative one. In other words, the bays and the leafs correspond to
the changes in the sign of the passage. Since our journey begins
with a positive passage and ends with a negative one, therefore
leafs $l_1,\dots,l_{p+1}$ and bays $b_1,\dots,b_p$ are visited in
the intertwining order $l_1,b_1,l_2,b_2,\dots,l_{p},b_{p},l_{p+1}$.
The number of the leafs (with the last leaf $l_{p+1}$ excluded) is
equal to the number of the bays, we can therefore consider a pairing
between them given by $l_i\mapsto b_i$ for $1\leq i\leq p$. In other
words, to a leaf $l$ we assign the first bay which is visited in our
journey after leaving $l$.

\subsection{Catalan sequences}
We say that $\epsilon=\big(\epsilon(1),\dots,\epsilon(k)\big)$ is a
Catalan sequence if $\epsilon(1),\dots,\epsilon(k)\in\{-1,+1\}$,
$\epsilon(1)+\cdots+\epsilon(k)=0$ and all partial sums are
non-negative: $\epsilon(1)+\cdots+\epsilon(l)\geq 0$ for all $1\leq
l\leq k$.

%

%
If $\epsilon$ is a Catalan sequence then there is no vertex $v\in
T_{\sigma}$ such that $v\prec R$.

\begin{lemma} \label{lem:catalanglue}
For  a Catalan sequence $\epsilon$ there exists a unique compatible
pairing $\sigma$ with the property that $R\preceq v$ for every
vertex $v\in T_{\sigma}$. We call it \emph{Catalan pairing}.

\end{lemma}
\begin{proof}
In the
sequence $\epsilon$ let us replace each element $+1$ by a left
bracket ``$\langle$" and let us replace each element $-1$ by a right
bracket ``$\rangle$". We leave it as an exercise to the Reader to
check that the pairing $\sigma$ between corresponding pairs of left
and right brackets is the unique pairing with the required property.
\end{proof}

\subsection{The main bijection}
The main result of this section is the algorithm 
\MainBijection{T} (with the auxiliary algorithm \SmallBijection{T}) which
provides the bijection announced in Theorem \ref{theo:mainresult}. In the
remaining part of the article we will show that this algorithm indeed provides
the desired bijection.

\begin{remark}
At the beginning of each iteration of the loop in \MainBijection $T$ is a
quotient tree $T_{\sigma}$ for some pairing $\sigma$ which is compatible with 
$\epsilon_i$. In order to check it (formally: by induction) we observe that
$l_i$ edges from each
side of the root in the polygonal graph $G_{\epsilon_i}$
are among those which were unglued in line \ref{linia:drugi-krok} of 
\MainBijection. These are the edges which we remove in
\ref{linia:remove-edges} of
\MainBijection. Formally, it corresponds to removal of the
first $l_i$ and the last $l_i$ elements from the
sequence $\epsilon_i$ and it is easily checked that the result is equal to
$(-\epsilon_{i-1})$. The change of the orientations of the edges in line 
\ref{linia:flip} means the change of sign of the corresponding sequence
$\epsilon$, hence after the iteration of the main loop in \MainBijection
$T=T_{\sigma}$ is a quotient tree corresponding to $\epsilon_{i-1}$.
\end{remark}

\begin{remark}
The operation of reversing the order $<$ in line \ref{linia:flip} of
\MainBijection means that we do not change the labels assigned to the tree $T$
but we change (by reversing) the way we compare them. It follows 
that for in line \ref{linia:malabijekcja} and in the function
\SmallBijection{T} we consider the set of labels (which is the set of
integer numbers) with its
usual order $<$ if $m-i$ is even and with the reverse of its
usual order if $m-i$ is odd.
\end{remark}

\begin{remark}
Tree $T$ in the algorithms \MainBijection and \SmallBijection
is always a quotient tree $T_{\sigma}$ for some pairing $\sigma$ which is
compatible with some sequence $\epsilon$. Each edge of this tree was created
from a pair of the edges of the polygonal graph $G_{\epsilon}$; therefore
the operation of ungluing in line \ref{ln:unglue} of \SmallBijection should be
understood as ungluing of these original edges. On a formal level ungluing
corresponds to removal of some pairs from the pairing $\sigma$; similarly
regluing in line \ref{ln:glue} of \SmallBijection corresponds to a
creation of new pairs in $\sigma$. Similar remarks concern lines
\ref{linia:drugi-krok}, \ref{lin:catalanglue} of
\MainBijection. 
\end{remark}

\begin{remark}
Lines \ref{linia:liczezatoke}--\ref{linia:liczezatoke2} of \SmallBijection
compute the bay $BA$, $CA$ corresponding in the tree $U$ to the leaf $D$.
\end{remark}


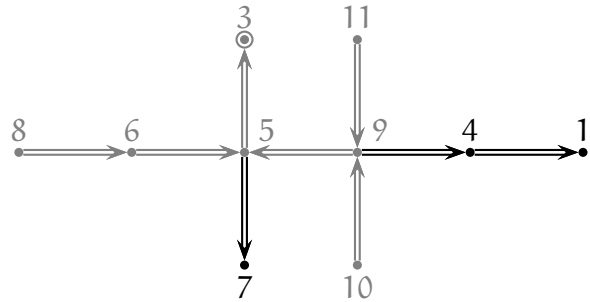
\begin{figure}
\psset{unit=1.5cm}
\begin{pspicture}[](-0.1,-1.3)(5.1,1.3)
 \cnode*[linecolor=gray](0,0){0.6mm}{A}
 \rput[u](0,0.2){$\gray 8$}

 \cnode*[linecolor=gray](1,0){0.6mm}{B}
 \rput[u](1,0.2){$\gray 6$}

 \cnode*[linecolor=gray](2,0){0.6mm}{C}
 \rput[u](2.2,0.2){$\gray 5$}

 \cnode*[linecolor=gray](3,0){0.6mm}{D}
 \rput[u](3.2,0.2){$\gray 9$}

 \cnode*(4,0){0.6mm}{E}
 \rput[u](4.0,0.2){$4$}

 \cnode*(5,0){0.6mm}{F}
 \rput[u](5,0.2){$1$}

 \cnode*[linecolor=gray](2,1){0.6mm}{G}
 \rput[u](2,1.2){$\gray 3$}
 \cnode[linecolor=gray](2,1){1.2mm}{G}

 \cnode*[linecolor=gray](3,1){0.6mm}{H}
 \rput[u](3,1.2){$\gray 11$}

 \cnode*(2,-1){0.6mm}{I}
 \rput[u](2,-1.2){$7$}

 \cnode*[linecolor=gray](3,-1){0.6mm}{J}
 \rput[u](3,-1.2){$\gray 10$}

 \ncline[arrowsize=2mm,doubleline=true,linecolor=gray]{->}{A}{B}
 \ncline[arrowsize=2mm,doubleline=true,linecolor=gray]{->}{B}{C}
 \ncline[arrowsize=2mm,doubleline=true]{->}{D}{E}
 \ncline[arrowsize=2mm,doubleline=true]{->}{E}{F}
 \ncline[arrowsize=2mm,doubleline=true,linecolor=gray]{->}{C}{G}
 \ncline[arrowsize=2mm,doubleline=true]{->}{C}{I}

 \ncline[arrowsize=2mm,doubleline=true,linecolor=gray]{<-}{C}{D}
 \ncline[arrowsize=2mm,doubleline=true,linecolor=gray]{<-}{D}{H}
 \ncline[arrowsize=2mm,doubleline=true,linecolor=gray]{<-}{D}{J}

 \end{pspicture}
 \caption[]{Algorithm \MainBijection{T}, line \ref{linia:pierwszy-krok}.
Subtree $U$ was marked in gray.}
 \label{fig:exampleTsigmaprim}
\end{figure}

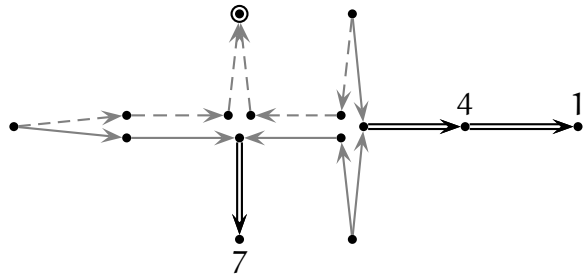
\begin{figure}
\psset{unit=1.5cm}
\begin{pspicture}[](-0.1,-1.3)(5.1,1.1)
 \cnode*(0,0){0.6mm}{A}
 \cnode*(1,0.1){0.6mm}{BN}
 \cnode*(1,-0.1){0.6mm}{BS}

 \cnode*(1.9,0.1){0.6mm}{CNW}
 \cnode*(2,-0.1){0.6mm}{CS}
 \cnode*(2.1,0.1){0.6mm}{CNE}

 \cnode*(2.9,0.1){0.6mm}{DNW}
 \cnode*(2.9,-0.1){0.6mm}{DSW}
 \cnode*(3.1,0){0.6mm}{DE}

 \cnode*(4,0){0.6mm}{E}

 \cnode*(5,0){0.6mm}{F}
 \cnode*(2,1){0.6mm}{G}
 \cnode(2,1){1.2mm}{G}
 \cnode*(3,1){0.6mm}{H}
 \cnode*(2,-1){0.6mm}{I}
 \cnode*(3,-1){0.6mm}{J}

\rput[u](4.0,0.2){$4$}
\rput[u](2,-1.2){$7$}
\rput[u](5,0.2){$1$}

 \ncline[arrowsize=2mm,linestyle=dashed,linecolor=gray]{->}{A}{BN}
 \ncline[arrowsize=2mm,linecolor=gray]{->}{A}{BS}
 \ncline[arrowsize=2mm,linestyle=dashed,linecolor=gray]{->}{BN}{CNW}
 \ncline[arrowsize=2mm,linecolor=gray]{->}{BS}{CS}
 \ncline[arrowsize=2mm,doubleline=true]{->}{CS}{I}
 \ncline[arrowsize=2mm,doubleline=true]{->}{DE}{E}
 \ncline[arrowsize=2mm,doubleline=true]{->}{E}{F}
 \ncline[arrowsize=2mm,linestyle=dashed,linecolor=gray]{->}{CNE}{G}
 \ncline[arrowsize=2mm,linestyle=dashed,linecolor=gray]{->}{CNW}{G}

 \ncline[arrowsize=2mm,linestyle=dashed,linecolor=gray]{<-}{CNE}{DNW}
 \ncline[arrowsize=2mm,linecolor=gray]{<-}{CS}{DSW}
 \ncline[arrowsize=2mm,linestyle=dashed,linecolor=gray]{<-}{DNW}{H}
 \ncline[arrowsize=2mm,linecolor=gray]{<-}{DE}{H}
 \ncline[arrowsize=2mm,linecolor=gray]{<-}{DSW}{J}
 \ncline[arrowsize=2mm,linecolor=gray]{<-}{DE}{J}

 \end{pspicture}
 \caption[]{Algorithm \MainBijection{T}, line \ref{linia:drugi-krok}.}
 \label{fig:Tprim}
\end{figure}

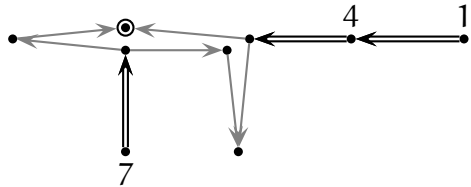
\begin{figure}
\psset{unit=1.5cm}
\begin{pspicture}[](-0.1,-1.3)(5.1,0.3)
 \cnode*(1,0){0.6mm}{B}

  \cnode*(2,0.1){0.6mm}{CN}
 \cnode(2,0.1){1.2mm}{CN}


 \cnode*(2,-0.1){0.6mm}{CS}

 \cnode*(2.9,-0.1){0.6mm}{DSW}
 \cnode*(3.1,0){0.6mm}{DE}

 \cnode*(4,0){0.6mm}{E}

 \cnode*(5,0){0.6mm}{F}
 \cnode*(2,-1){0.6mm}{I}
 \cnode*(3,-1){0.6mm}{J}

 \ncline[arrowsize=2mm,linecolor=gray]{->}{B}{CN}
 \ncline[arrowsize=2mm,linecolor=gray]{->}{CS}{DSW}
 \ncline[arrowsize=2mm,linecolor=gray]{->}{DSW}{J}
 \ncline[arrowsize=2mm,linecolor=gray]{->}{DE}{J}

 \ncline[arrowsize=2mm,linecolor=gray]{<-}{CN}{DE}
 \ncline[arrowsize=2mm,doubleline=true]{<-}{CS}{I}
 \ncline[arrowsize=2mm,doubleline=true]{<-}{DE}{E}
 \ncline[arrowsize=2mm,doubleline=true]{<-}{E}{F}
 \ncline[arrowsize=2mm,linecolor=gray]{<-}{B}{CS}

\rput[u](4.0,0.2){$4$}
\rput[u](2,-1.2){$7$}
\rput[u](5,0.2){$1$}


 \end{pspicture}

 \caption[]{Algorithm \MainBijection{T}, line \ref{linia:trzeci-krok}. This
graph was obtained from  Figure
\ref{fig:Tprim} by removal of  the dashed edges and it can be regarded
as a certain polygonal graph $G_{\epsilon'}$ with a number of trees attached to
it. The orientation of all edges was reversed.}

 \label{fig:Tbis}
\end{figure}

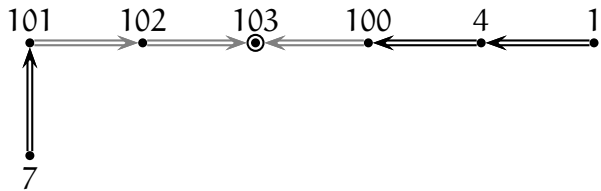
\begin{figure}
\psset{unit=1.5cm}
\begin{pspicture}[](-0.1,-1.1)(5.1,0.4)

 \cnode*(0,-1){0.6mm}{Z}
 \rput[u](0,-1.2){$7$}

 \cnode*(0,0){0.6mm}{A}
 \rput[u](0,0.2){$101$}

 \cnode*(1,0){0.6mm}{B}
 \rput[u](1,0.2){$102$}

 \cnode*(2,0){0.6mm}{C}
 \cnode(2,0){1.2mm}{C}
 \rput[u](2,0.2){$103$}

 \cnode*(3,0){0.6mm}{D}
 \rput[u](3,0.2){$100$}

 \cnode*(4,0){0.6mm}{E}
 \rput[u](4,0.2){$4$}

 \cnode*(5,0){0.6mm}{F}
 \rput[u](5,0.2){$1$}

 \ncline[arrowsize=2mm,doubleline=true,linecolor=gray]{->}{A}{B}
 \ncline[arrowsize=2mm,doubleline=true,linecolor=gray]{->}{B}{C}

 \ncline[arrowsize=2mm,doubleline=true]{<-}{A}{Z}
 \ncline[arrowsize=2mm,doubleline=true,linecolor=gray]{<-}{C}{D}
 \ncline[arrowsize=2mm,doubleline=true]{<-}{D}{E}
 \ncline[arrowsize=2mm,doubleline=true]{<-}{E}{F}

 \end{pspicture}

 \caption[]{Algorithm \MainBijection{T}, line \ref{linia:koniec}. The
polygonal graph $G_{\epsilon'}$ from figure \ref{fig:Tbis} was glued
according to the Catalan pairing. Artificial labels $100$--$103$ were created to
label new vertices. The order of the labels was reversed therefore
$103<102<101<100<7<4<1$.}

  \label{fig:Tbabuga}
\end{figure}

\begin{function}
\SetLine
\SetKwInOut{Input}{input}
\SetKwInOut{Output}{output}
\Input{$l_1\leq \cdots \leq l_m$ positive integers, $L=l_1+\cdots+l_m+1$ \\
$T$ is a quotient tree corresponding to the Catalan sequence $\epsilon_m$ \\
equipped with a total  order $<$ which is compatible with the \\ orientations of
edges}
\Output{disjoint sets $B_1,\dots,B_m$ such that $B_1 \cup \cdots \cup
B_m=\{1,\dots,L\}$ and
$$ |B_1|+\cdots+|B_n|\leq l_1+\cdots+l_n$$
holds true for each $1\leq n \leq m-1$}
\SetKw{downto}{downto}
\SetKwFunction{Sink}{SmallBijection}
\SetKwData{Mylabels}{labels}

\BlankLine

label all vertices of $T$ with numbers $1,\dots,L$ in such a way that each
label appears exactly once and the order
$<$ of vertices coincides with the order of the labels\;

\For{i=m \downto 1}{

$T\leftarrow\Sink(T)$\; \nllabel{linia:malabijekcja}
$U\leftarrow \text{tree } \{x\in T: x\succeq R\}$\;
$B_i\leftarrow(\text{labels of the vertices of }U)\cap \{1,\dots,L\} $\; 
\nllabel{linia:pierwszy-krok}
\tcc{cf.\ figure \ref{fig:exampleTsigmaprim}}
\BlankLine

remove the labels of the vertices of $U$\; 
\nllabel{lin:remove}   
unglue all edges of tree $U$\;
\nllabel{linia:drugi-krok}
\tcc{cf.\ figure \ref{fig:Tprim}}

\BlankLine

remove $l_i$ edges at each side of the vertex $R$\; \nllabel{linia:remove-edges}
change the orientation of all edges and reverse the order $<$\;
\nllabel{linia:flip}
\tcc{cf.\ figure \ref{fig:Tbis}}
\nllabel{linia:trzeci-krok}

\BlankLine

create sufficiently many artifical labels (integer numbers all different
from $1,\dots,L$) which are
smaller than any label on tree $T$\;

glue the remaining edges of tree $U$ by the Catalan pairing given by Lemma
\ref{lem:catalanglue}\; \nllabel{lin:catalanglue}
label the unlabeled vertices with artificial labels in such a way that on tree
$U$ the orders $<$ and $\lhd$ coincide\; \nllabel{linia:koniec}
\tcc{cf.\ figure \ref{fig:Tbabuga}}

}
\Return $B_1,\dots,B_m$\;

%
\caption{MainBijection(T)}
\label{MainBijection}
\end{function}

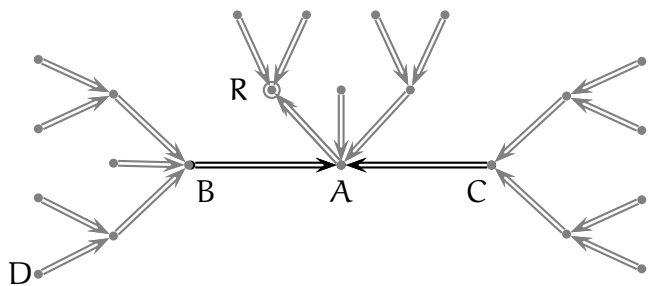
\begin{figure}
 \psset{unit=2cm}
 \begin{pspicture}[](-1,-0.8)(1,1.2)
 \cnode*(0,0){0.6mm}{A}
 \cnode*(1,0){0.6mm}{C}
 \cnode*(-1,0){0.6mm}{B}
 \ncline[arrowsize=2mm,doubleline=true]{<-}{A}{B}
 \ncline[arrowsize=2mm,doubleline=true]{<-}{A}{C}
 \rput[t](0,-0.1){$A$}
 \rput[t](0.9,-0.1){$C$}
 \rput[t](-0.9,-0.1){$B$}

\psset{linecolor=gray}

\renewcommand{\psedge}[2]{\ncline[doubleline=true,arrowsize=2mm]{->}{#1}{#2}}
 \rput[B](0,0){\pstree[treemode=U,levelsep=0.5,treesep=0.4]{\Tc*{0.6mm}}{
     \pstree{ \Tc{1.2mm}~[tnpos=l,tnsep=2mm]{$R$} \Tc*{0.6mm}  }{
    
\renewcommand{\psedge}[2]{\ncline[doubleline=true,arrowsize=2mm]{<-}{#1}{#2}}
   \Tc*{0.6mm}  \Tc*{0.6mm}}
   \renewcommand{\psedge}[2]{\ncline[doubleline=true,arrowsize=2mm]{<-}{#1}{#2}}
   \Tc*{0.6mm}\pstree{\Tc*{0.6mm}}{\Tc*{0.6mm} \Tc*{0.6mm}}}}

 \renewcommand{\psedge}[2]{\ncline[doubleline=true,arrowsize=2mm]{<-}{#1}{#2}}
 \rput[B](-1.6,0){\pstree[treemode=L,levelsep=0.5,treesep=0.4]{\Tc*{0.6mm}}{
     \pstree{\Tc*{0.6mm}}{\Tc*{0.6mm} \Tc*{0.6mm}}
\Tc*{0.6mm}\pstree{\Tc*{0.6mm}}{\Tc*{0.6mm}
     \Tc*{0.6mm}~[tnpos=l,tnsep=0mm]{$D$}    }}}

 \rput[B](1.5,0){\pstree[treemode=R,levelsep=0.5,treesep=0.4]{\Tc*{0.6mm}}{
     \pstree{\Tc*{0.6mm}}{\Tc*{0.6mm}
     \Tc*{0.6mm}} \pstree{\Tc*{0.6mm}}{\Tc*{0.6mm} \Tc*{0.6mm} }}}

 \end{pspicture}
\SetKwFunction{SmallBijection}{SmallBijection}
\caption[]{Algorithm \SmallBijection{T}, case $D\neq B$. The order of the
vertices is given by $R\leq A<B<C<D$. Note that only edges belonging to the
subtree $U$ are displayed.}

 \label{fig:transform0}
\end{figure}

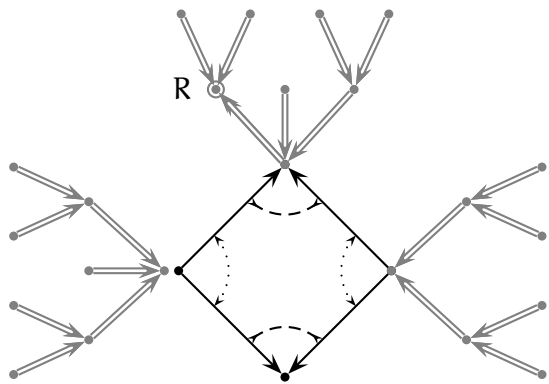
\begin{figure}
\psset{unit=2cm}
\begin{pspicture}[](-1,-0.8)(1,1.8)
 \cnode*(0.7071067,0){0.6mm}{E}
 \cnode*(0,0.7071067){0.6mm}{N}
 \cnode*(-0.7071067,0){0.6mm}{W}
 \cnode*(0,-0.7071067){0.6mm}{S}
 \psarc[linestyle=dashed]{>-<}(0,0.7071067){0.33333}{225}{315}
 \psarc[linestyle=dashed]{>-<}(0,-0.7071067){0.33333}{45}{135}
 \psarc[linestyle=dotted]{<->}(0.7071067,0){0.33333}{135}{225}
 \psarc[linestyle=dotted]{<->}(-0.7071067,0){0.33333}{-45}{45}
 \ncline[arrowsize=2mm,doubleline=false]{<-}{N}{E}
 \ncline[arrowsize=2mm,doubleline=false]{<-}{N}{W}
 \ncline[arrowsize=2mm,doubleline=false]{<-}{S}{E}
 \ncline[arrowsize=2mm,doubleline=false]{<-}{S}{W}

\psset{linecolor=gray}

\renewcommand{\psedge}[2]{\ncline[doubleline=true,arrowsize=2mm]{->}{#1}{#2}}
 \rput[B](0,0.7071067){\pstree[treemode=U,levelsep=0.5,treesep=0.4]{\Tc*{0.6mm}}
{
     \pstree{ \Tc{1.2mm}~[tnpos=l,tnsep=2mm]{$R$} \Tc*{0.6mm}  }{
    
\renewcommand{\psedge}[2]{\ncline[doubleline=true,arrowsize=2mm]{<-}{#1}{#2}}
   \Tc*{0.6mm}  \Tc*{0.6mm}}
   \renewcommand{\psedge}[2]{\ncline[doubleline=true,arrowsize=2mm]{<-}{#1}{#2}}
   \Tc*{0.6mm}\pstree{\Tc*{0.6mm}}{\Tc*{0.6mm} \Tc*{0.6mm}}}}

 \renewcommand{\psedge}[2]{\ncline[doubleline=true,arrowsize=2mm]{<-}{#1}{#2}}
 \rput[B](-1.303,0){\pstree[treemode=L,levelsep=0.5,treesep=0.4]{\Tc*{0.6mm}}{
     \pstree{\Tc*{0.6mm}}{\Tc*{0.6mm} \Tc*{0.6mm}}
\Tc*{0.6mm}\pstree{\Tc*{0.6mm}}{\Tc*{0.6mm}
     \Tc*{0.6mm}~[tnpos=l,tnsep=0mm]{
}    }}}

 \rput[B](+1.2071067,0){\pstree[treemode=R,levelsep=0.5,treesep=0.4]{\Tc*{0.6mm}
}{
     \pstree{\Tc*{0.6mm}}{\Tc*{0.6mm}
     \Tc*{0.6mm}} \pstree{\Tc*{0.6mm}}{\Tc*{0.6mm} \Tc*{0.6mm}}}}
 \end{pspicture}
 \caption{The tree from Figure \ref{fig:transform0} after ungluing the edges
$BA$ and $CA$.}
 \label{fig:transform1}
\end{figure}

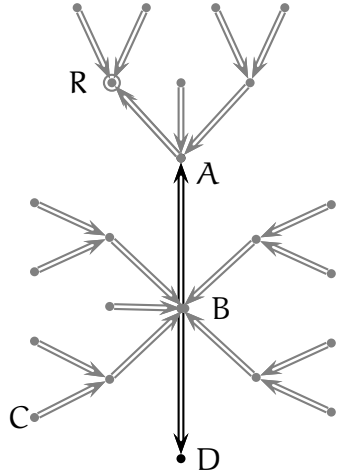
\begin{figure}
\psset{unit=2cm}
\begin{pspicture}[](-1,-1.1)(1,2)
 \cnode*(0,1){0.6mm}{A}
 \cnode*(0,0){0.6mm}{B}
 \cnode*(0,-1){0.6mm}{S}

 \ncline[arrowsize=2mm,doubleline=true]{<-}{A}{B}
 \ncline[arrowsize=2mm,doubleline=true
]{<-}{S}{B}
 \rput[l](0.1,0.9){$A$}
 \rput[l](0.2,0){$B$}
 \rput[l](0.1,-1){$D$}

\psset{linecolor=gray}
\renewcommand{\psedge}[2]{\ncline[doubleline=true,arrowsize=2mm]{->}{#1}{#2}}
 \rput[B](0,1){\pstree[treemode=U,levelsep=0.5,treesep=0.4]{\Tc*{0.6mm}}{
     \pstree{ \Tc{1.2mm}~[tnpos=l,tnsep=2mm]{$R$} \Tc*{0.6mm}  }{
    
\renewcommand{\psedge}[2]{\ncline[doubleline=true,arrowsize=2mm]{<-}{#1}{#2}}
   \Tc*{0.6mm}  \Tc*{0.6mm}}
   \renewcommand{\psedge}[2]{\ncline[doubleline=true,arrowsize=2mm]{<-}{#1}{#2}}
   \Tc*{0.6mm}\pstree{\Tc*{0.6mm}}{\Tc*{0.6mm} \Tc*{0.6mm}}}}

 \renewcommand{\psedge}[2]{\ncline[doubleline=true,arrowsize=2mm]{<-}{#1}{#2}}
 \rput[B](-0.55,0){\pstree[treemode=L,levelsep=0.5,treesep=0.4]{\Tc*{0.6mm}}{
     \pstree{\Tc*{0.6mm}}{\Tc*{0.6mm} \Tc*{0.6mm}}
\Tc*{0.6mm}\pstree{\Tc*{0.6mm}}{\Tc*{0.6mm}
     \Tc*{0.6mm}~[tnpos=l,tnsep=0mm]{$C$}    }}}

 \rput[B](0.5,0){\pstree[treemode=R,levelsep=0.5,treesep=0.4]{\Tc*{0.6mm}}{
     \pstree{\Tc*{0.6mm}}{\Tc*{0.6mm}
     \Tc*{0.6mm}}  \pstree{\Tc*{0.6mm}}{\Tc*{0.6mm} \Tc*{0.6mm}}}}
 \end{pspicture}
 \caption{The tree from Figure \ref{fig:transform0} after regluing the edges
$BA$ and
$CA$ in a different way. Please notice the change of the labels of
the vertices $A,B,C,D$. }

 \label{fig:transform2}
\end{figure}

\begin{figure}
 \psset{unit=2cm}
 \begin{pspicture}[](-1,-0.8)(1,1.2)
 \cnode*(0,0){0.6mm}{A}
 \cnode*(1,0){0.6mm}{C}
 \cnode*(-1,0){0.6mm}{B}
 \ncline[arrowsize=2mm,doubleline=true]{<-}{A}{B}
 \ncline[arrowsize=2mm,doubleline=true]{<-}{A}{C}
 \rput[t](0,-0.1){$A$}
 \rput[t](0.9,-0.1){$C$}
 \rput[t](-1,-0.1){$D$}

\psset{linecolor=gray}

\renewcommand{\psedge}[2]{\ncline[doubleline=true,arrowsize=2mm]{->}{#1}{#2}}
 \rput[B](0,0){\pstree[treemode=U,levelsep=0.5,treesep=0.4]{\Tc*{0.6mm}}{
     \pstree{ \Tc{1.2mm}~[tnpos=l,tnsep=2mm]{$R$} \Tc*{0.6mm}  }{
\renewcommand{\psedge}[2]{\ncline[doubleline=true,arrowsize=2mm]{<-}{#1}{#2}}
   \Tc*{0.6mm}  \Tc*{0.6mm}}
   \renewcommand{\psedge}[2]{\ncline[doubleline=true,arrowsize=2mm]{<-}{#1}{#2}}
   \Tc*{0.6mm}\pstree{\Tc*{0.6mm}}{\Tc*{0.6mm} \Tc*{0.6mm}}}}

  \renewcommand{\psedge}[2]{\ncline[doubleline=true,arrowsize=2mm]{<-}{#1}{#2}}
 \rput[B](1.5,0){\pstree[treemode=R,levelsep=0.5,treesep=0.4]{\Tc*{0.6mm}}{
     \pstree{\Tc*{0.6mm}}{\Tc*{0.6mm}
     \Tc*{0.6mm}} \pstree{\Tc*{0.6mm}}{\Tc*{0.6mm} \Tc*{0.6mm}}}}

 \end{pspicture}
\caption[]{Algorithm \SmallBijection{T}, case $D=B$. The order of vertices is
given by $R\leq A<C<D$.}
 \label{fig:transform0bis}
\end{figure}
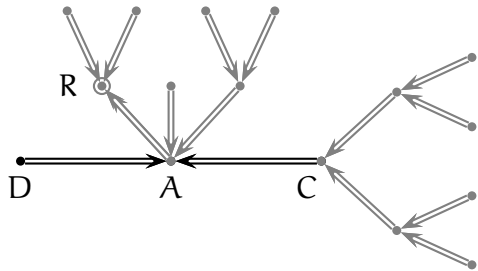

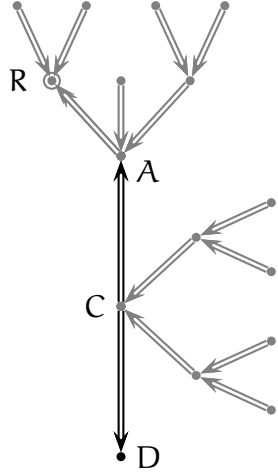
\begin{figure}
\psset{unit=2cm}
\begin{pspicture}[](-1,-1.1)(1,2)
 \cnode*(0,1){0.6mm}{A}
 \cnode*(0,0){0.6mm}{B}
 \cnode*(0,-1){0.6mm}{S}

 \ncline[arrowsize=2mm,doubleline=true]{<-}{A}{B}
 \ncline[arrowsize=2mm,doubleline=true
]{<-}{S}{B}
 \rput[l](0.1,0.9){$A$}
 \rput[r](-0.1,0){$C$}
 \rput[l](0.1,-1){$D$}

\psset{linecolor=gray}
\renewcommand{\psedge}[2]{\ncline[doubleline=true,arrowsize=2mm]{->}{#1}{#2}}
 \rput[B](0,1){\pstree[treemode=U,levelsep=0.5,treesep=0.4]{\Tc*{0.6mm}}{
     \pstree{ \Tc{1.2mm}~[tnpos=l,tnsep=2mm]{$R$} \Tc*{0.6mm}  }{
    
\renewcommand{\psedge}[2]{\ncline[doubleline=true,arrowsize=2mm]{<-}{#1}{#2}}
   \Tc*{0.6mm}  \Tc*{0.6mm}}
   \renewcommand{\psedge}[2]{\ncline[doubleline=true,arrowsize=2mm]{<-}{#1}{#2}}
   \Tc*{0.6mm}\pstree{\Tc*{0.6mm}}{\Tc*{0.6mm} \Tc*{0.6mm}}}}

  \renewcommand{\psedge}[2]{\ncline[doubleline=true,arrowsize=2mm]{<-}{#1}{#2}}
 \rput[B](0.5,0){\pstree[treemode=R,levelsep=0.5,treesep=0.4]{\Tc*{0.6mm}}{
     \pstree{\Tc*{0.6mm}}{\Tc*{0.6mm}
     \Tc*{0.6mm}}  \pstree{\Tc*{0.6mm}}{\Tc*{0.6mm} \Tc*{0.6mm}}}}
 \end{pspicture}
 \caption{The tree from Figure \ref{fig:transform0bis} after regluing the edges
$DA$ and $CA$
 in a different way. Please notice the change of the labels of the vertices
$A,C,D$.}
 \label{fig:transform2bis}
\end{figure}

\begin{function}
\SetLine
\SetKwInOut{Input}{input}
\SetKwInOut{Output}{output}
\Input{$T$ is a quotient
tree correspondig to some sequence $\epsilon$. The vertices of $T$
are equipped with some labels in such a way
that the order of labels is compatible with the orientations of the edges.}
\Output{Tree $T$ which is a quotient tree corresponding to the same 
sequence $\epsilon$. The vertices of $T$ are labeled; the set of labels of this
output tree coincides with the set of the vertices of the input tree. More
detailed description of the output will be given in Theorem
\ref{theo:malabijekcja}.}
\SetKwFunction{father}{father}
\SetKwFunction{Sink}{Sink}
\SetKwData{Mylabels}{labels}

\BlankLine

%

\While{orders $<$ and $\lhd$ do not coincide on $\{x\in T: x\succeq R\}$}{
$D\leftarrow$the minimal element (with respect to $<$) such that $R\prec D$ and
orders $<$ and $\lhd$ do not coincide on $\{x\in T: R\preceq x \text{ and } x
\leq D\}$ \; \nllabel{lin:d-rowna-sie}

$U\leftarrow\text{tree }\{x\in T: R\preceq x \text{ and } x \leq D\}$ \;

$C\leftarrow$the successor of $D$ in $U$ with respect to $\lhd$\;
\nllabel{linia:liczezatoke}
$A\leftarrow$father of $C$\;
$B\leftarrow$son of $A$ in $U$ which is to the left of $C$\;
\nllabel{linia:liczezatoke2}
$\Mylabels\leftarrow$set of labels carried by the vertices $A$, $B$, $C$, $D$\;
\tcc{
cf.\ figure \ref{fig:transform0} and figure
\ref{fig:transform0bis}}

\BlankLine

remove the labels from the vertices $A$, $B$, $C$, $D$\;
unglue the edges $BA$ and $CA$ \nllabel{ln:unglue}
\tcc*{cf. figure  \ref{fig:transform1}}

\BlankLine

reglue these edges in the other possible way\; \nllabel{ln:glue}
to unlabeled vertices give labels from $\Mylabels$ in such a way that for each
pair of newly labeled vertices $x<y$ iff $x\lhd y$\; \nllabel{ln:koniec-malego}
\tcc{cf.\ figure  \ref{fig:transform2} and figure \ref{fig:transform2bis}}
   }

%

\Return{$T$}\;

\caption{SmallBijection(T)}
\end{function}

\section{Proof of the correctness of the small bijection}
\label{sec:almostthemain}

\subsection{Statement of the result}
\label{subsec:basicbijection}


Let $T_{\sigma}$ be a quotient tree and let $<$ be an order on its vertices. We
may always label the vertices of $T_{\sigma}$ (for example, with integer
numbers) in such a way that the order of the vertices coincides with the order
of the corresponding labels. In this way we can view $\SmallBijection$ as a map
which to a pair $(T_{\sigma},<)$ (or, more formally, $(\sigma,<)$) associates
another pair of this form. 

\begin{theorem}
\label{theo:malabijekcja} Let $\epsilon=\big( \epsilon(1),\dots,
\epsilon(k) \big)$ be a Catalan sequence.
The function $\SmallBijection$ as described above provides a
bijection between
\begin{enumerate}
\item[(A)]  the set of pairs $(\sigma,<)$, where $\sigma$ is a
pairing compatible with $\epsilon$ and $<$ is a total order on the
vertices of \/ $T_{\sigma}$ compatible with the orientation of the
edges;


\item[(B)]  the set of pairs $(\sigma,<)$, where $\sigma$ is a
pairing compatible with $\epsilon$ and $<$ is a total order on the
vertices of \/ $T_{\sigma}$ with the following two properties:
\label{bij2}

\begin{itemize}
\item  on the set $\{x\in T_{\sigma}: x\succeq R\}$ the orders $<$ and $\lhd$
coincide;

\item for all pairs of vertices $v,w\in
T_{\sigma}$ such that $R \not\preceq v$ and $R\not\preceq w$ we have
$$v \prec w \implies v< w.$$

\end{itemize}

%
%
\end{enumerate}
\end{theorem}

The remaining part of this section is devoted to the proof of this theorem.


\subsection{Intermediate triples}
\label{subsec:intermediatetriples} 
Our strategy is to describe precisely which pairs $(\sigma,<)$ (or
alternatively, trees $T$) might arise in the intermediate steps of algorithm
\SmallBijection.


\begin{definition}
We call $(\sigma,<,S)$ an intermediate triple if $\sigma$ is a
pairing compatible with $\epsilon$, $<$ is a total order on the
vertices of $T_{\sigma}$ and $S$ is one of the vertices of
$T_{\sigma}$ with the following properties:
\begin{enumerate}
\item $R\preceq S$ and $R\leq S$, where $R$ denotes the root; \label{interm1}

\item on the set $\{x: R\preceq x \text{ and } x \leq S\}$ the
orders $<$ and $\lhd$ coincide;
\label{interm2}

\item for all pairs of adjacent vertices $v,w\in T_{\sigma}$ such that $v\prec w$ and
$v>w$ we have $R \not\preceq v$ and $R\prec w$ and the set $\{x\in
T_{\sigma}: R\preceq x\text{ and } S<x<v\}$ is empty.
\label{interm3}

\end{enumerate}

\end{definition}

\subsection{Startpoints and endpoints}

\begin{lemma}
\label{lem:identyfikacjaA} Intermediate triples $(\sigma,<,S)$ for
which $S=R$ are in a one-to-one correspondence with the pairs
$(\sigma,<)$ which contribute to the set (A) and thus to the possible input
data of algorithm \SmallBijection.
\end{lemma}
\begin{proof}
Suppose that $(\sigma,<)$ contributes to (A); we set $S=R$. In order
to show property (\ref{interm2}) of intermediate triples it is enough to
observe that if $x$
fulfills $R\preceq x$ then also $R\leq x$, therefore the set $\{x\in
T_{\sigma}: R \preceq x \text{ and } x\leq R\}$ consists of only one
element $R$. The other two properties of intermediate triples hold
true trivially.

Suppose that $(\sigma,<,R)$ is an intermediate triple and suppose
that there exists a pair of vertices $v,w$ such that $v\prec w$ and
$v>w$. With no loss of generality we may assume that the vertices
$v$ and $w$ are adjacent (if this is not the case we may find a pair
of adjacent vertices $v'$, $w'$ such that $v'\prec w'$ and $v'>w'$
on the path from the vertex $v$ to the vertex $w$). In the case
$w\leq R$ property (\ref{interm3}) shows that $R\prec w$ and
property (\ref{interm2}) shows that since $R\lhd w$ therefore $R<w$
which contradicts $w\leq R$. In the case $w>R$ the set $\{x\in
T_{\sigma}: R\preceq x \text{ and } R<x<v\}$ contains $w$ which
contradicts property (\ref{interm3}). In this way we proved that the
total order $<$ is compatible with the orientations of the edges.

\end{proof}

\begin{lemma}
\label{lem:identyfikacjaB} Intermediate triples $(\sigma,<,S)$ for
which $S$ is the maximal element (with respect to the order $<$) of
the set $\{x\in T_{\sigma}: x\succeq R\}$ are in a one-to-one
correspondence with the pairs $(\sigma,<)$ which contribute to the
set (B). For such values of $T=T_{\sigma}$ algorithm \SmallBijection
terminates.
\end{lemma}
\begin{proof}
Suppose that $(\sigma,<)$ contributes to the set (B); we define $S$
to be the maximal element (with respect to $<$) of the set $\{x\in
T_{\sigma}: x\succeq R\}$. In order to prove that $(\sigma,<,S)$ is
an intermediate triple it is enough to show property
(\ref{interm3}) of intermediate triples. Suppose that there exist adjacent
vertices $v,w$,
such that $v\prec w$ and $v>w$; we shall consider now three cases.
The first case, $R\not\preceq v,w$ is not possible, since then $v<w$
would contradict $v>w$. The second case, $R\preceq v\prec w$ would
imply $v\lhd w$ and hence $v<w$ again contradicts $v>w$. Therefore,
the only remaining possibility is $R\preceq w$ and $R\not\preceq v$.
It is not possible that $R=w$ since then $v\prec R$ would contradict
the assumption that $\epsilon$ is a Catalan sequence. In this way we
proved that $R\prec w$, $R\not\preceq v$ which finishes the proof.

Suppose that $(\sigma,<,S)$ is as in the statement of the lemma. In
order to prove that $(\sigma,<)$ contributes to (B) it is enough to
prove that for all vertices $v,w$ such that $R\not\preceq v$ and
$R\not\preceq w$ we have $v\prec w \implies v<w$. If this is not the
case then there exist vertices $v,w$ such that $R\not\preceq w$ and
$v\prec w$, $v>w$. With no loss of generality we may assume that the
vertices $v$, $w$ are adjacent hence $R\prec w$ contradicts
$R\not\preceq w$.
\end{proof}

\subsection{The forward transformation}
\label{subsec:forward}

In this section we shall describe a certain invertible operation on
intermediate triples which will turn out to be equivalent to the algorithm
\SmallBijection. After a sufficient number of iterations every
intermediate triple corresponding to some element of (A) gets
transformed into an intermediate triple corresponding to some
element of (B). In this way the operation described in this section (which
coincides with \SmallBijection)
provides a bijection from (A) to (B).

Let $(\sigma,<,S)$ be an intermediate triple and let $D$ be the
smallest element (with respect to the order $<$) in the set $\{x:
x\succ R \text{ and } x>S\}$. If no such element $D$ exists, this
means that the triple $(\sigma,<,S)$ is as in Lemma
\ref{lem:identyfikacjaB} and hence can be identified with an element
of the set (B); in other words our algorithm finished its work.

If $(\sigma,<,D)$ is an intermediate triple, then we iterate our
procedure.

We consider now the opposite case when $(\sigma,<,D)$ is not an
intermediate triple. Then $S$ is the maximal element (with respect to $<$)
such that $R\preceq S$ and such that on the set $\{x\in T_{\sigma}: R\preceq
x\text{ and }
x\leq S\}$ the orders $<$ and $\lhd$ coincide. Also, $D$ is as prescribed by 
line \ref{lin:d-rowna-sie} of \SmallBijection.

Let us denote by $(\sigma',<)$ the pairing and the order which correspond to
the value of $T$ in line \ref{ln:koniec-malego} of algorithm
\SmallBijection.

\begin{lemma}
The triple $(\sigma',<,S)$ given by the above construction is an
intermediate triple.

The above procedure will stop after a finite number of steps.
\end{lemma}
\begin{proof}
We shall consider only the case when $B\neq D$, since the other case
is analogous. Conditions (\ref{interm1}) and (\ref{interm2}) are
very easy to verify. To check condition (\ref{interm3}) we need to
find adjacent pairs of vertices $v,w$ on tree $T_{\sigma'}$ for
which $v\prec w$ and $v>w$. Since condition (\ref{interm3}) is
fulfilled for the tree $T_{\sigma}$ it is enough to restrict our
attention to such pairs which are new, i.e.\ which were not present
on the tree $T_{\sigma}$. Figure \ref{fig:transform2} indicates one
such pair (namely $v=D$, $w=B$ and it is easy to check that this
pair causes no problems) there might be however some other such
pairs which were not shown on Figure \ref{fig:transform1} because
$v\notin U$ or $w\notin U$. It is easy to check that such new pairs
must fall into one of the following three categories: $C<v<D$, $w=C$
(causes no problems); or $v=D$, $A<w<D$ (impossible, since it would
imply that in the tree $T_{\sigma}$ the vertex $A$ is adjacent to
the vertex $w$ such that $A<w<D$ but $(\sigma,<,S)$ is an
intermediate triple, contradiction); or $B<v<C$, $w=B$ (causes no
problems).

Note that in each step of our operation the cardinality of $\{x\in
T_{\sigma}: R\preceq x \text{ and } S\leq x\}$ decreases. This shows
that our procedure will eventually stop.
\end{proof}

\subsection{The backward transformation}
\label{subsec:backward} In this section we shall describe the
inverse of the transformation from Section \ref{subsec:forward}.

Let $(\sigma,<,S)$ be an intermediate triple and let $S'$ be the
biggest element (with respect to the order $<$) of the set
$\{x:x\succeq R \text{ and } x<S\}$. If no such element exists it
means that $(\sigma,<,S)$ is as in Lemma \ref{lem:identyfikacjaA}
hence our algorithm finished its work. If $(\sigma,<,S')$ is an
intermediate triple, we can iterate our procedure.

We consider now the opposite case when $(\sigma,<,S')$ is not an
intermediate triple. It is possible only when condition
(\ref{interm3}) of an intermediate triple is not fulfilled, namely
there is a pair of adjacent vertices $B,D$ such that $B\succ D$ and
$B<D$ and $\{x: R\preceq x\text{ and } S'<x<D\}$ is a non-empty set.
There may be many pairs $B,D$ with this property; let us select the
one for which $D$ takes its maximal value (with respect to $<$).
Since $(\sigma,<,S)$ is an intermediate triple therefore $\{x:
R\preceq x\text{ and } S<x<D\}$ is empty. It follows that the only
element which could possibly belong to $\{x: R\preceq x\text{ and }
S'<x<D\}$ is equal to $S$ and therefore $S<D$. In particular,
\begin{equation} \{x: R\preceq x\text{ and } S<x<D\}=\emptyset.
\label{equ:nicmiedzyTiD}
\end{equation}

\begin{figure}
\psset{unit=2cm}
\begin{pspicture}[](-1,-1.1)(1,2)
 \cnode*(0,1){0.6mm}{A}
 \cnode*(0,0){0.6mm}{B}
 \cnode*(0,-1){0.6mm}{S}

 \ncline[arrowsize=2mm,doubleline=true]{<-}{A}{B}
 \ncline[arrowsize=2mm,doubleline=true,linestyle=dashed]{<-}{S}{B}
 \rput[l](0.1,0.9){$A$}
 \rput[r](-0.1,0){$B$}
 \rput[l](0.1,-1){$D$}

\psset{linecolor=gray}
\renewcommand{\psedge}[2]{\ncline[doubleline=true,arrowsize=2mm]{->}{#1}{#2}}
 \rput[B](0,1){\pstree[treemode=U,levelsep=0.5,treesep=0.4]{\Tc*{0.6mm}}{
     \pstree{ \Tc{1.2mm}~[tnpos=l,tnsep=2mm]{$R$} \Tc*{0.6mm}  }{
     \renewcommand{\psedge}[2]{\ncline[doubleline=true,arrowsize=2mm]{<-}{#1}{#2}}
   \Tc*{0.6mm}  \Tc*{0.6mm}}
   \renewcommand{\psedge}[2]{\ncline[doubleline=true,arrowsize=2mm]{<-}{#1}{#2}}
   \Tc*{0.6mm}\pstree{\Tc*{0.6mm}}{\Tc*{0.6mm} \Tc*{0.6mm}}}}

  \renewcommand{\psedge}[2]{\ncline[doubleline=true,arrowsize=2mm]{<-}{#1}{#2}}
 \rput[B](0.5,0){\pstree[treemode=R,levelsep=0.5,treesep=0.4]{\Tc*{0.6mm}}{
     \pstree{\Tc*{0.6mm}}{\Tc*{0.6mm}
     \Tc*{0.6mm}}  \pstree{\Tc*{0.6mm}}{\Tc*{0.6mm} \Tc*{0.6mm}}}}
 \end{pspicture}
 \hspace{2cm}
 \begin{pspicture}[](-1,-1.1)(1,2)
 \cnode*(0,1){0.6mm}{A}
 \cnode*(0,0){0.6mm}{B}
 \cnode*(0,-1){0.6mm}{S}

 \ncline[arrowsize=2mm,doubleline=true]{<-}{A}{B}
 \ncline[arrowsize=2mm,doubleline=true,linestyle=dashed]{->}{S}{B}
 \rput[l](0.1,0.9){$A$}
 \rput[r](-0.1,0){$B$}
 \rput[l](0.1,-1){$D$}

\psset{linecolor=gray}
\renewcommand{\psedge}[2]{\ncline[doubleline=true,arrowsize=2mm]{->}{#1}{#2}}
 \rput[B](0,1){\pstree[treemode=U,levelsep=0.5,treesep=0.4]{\Tc*{0.6mm}}{
     \pstree{ \Tc{1.2mm}~[tnpos=l,tnsep=2mm]{$R$} \Tc*{0.6mm}  }{
     \renewcommand{\psedge}[2]{\ncline[doubleline=true,arrowsize=2mm]{<-}{#1}{#2}}
   \Tc*{0.6mm}  \Tc*{0.6mm}}
   \renewcommand{\psedge}[2]{\ncline[doubleline=true,arrowsize=2mm]{<-}{#1}{#2}}
   \Tc*{0.6mm}\pstree{\Tc*{0.6mm}}{\Tc*{0.6mm} \Tc*{0.6mm}}}}

  \renewcommand{\psedge}[2]{\ncline[doubleline=true,arrowsize=2mm]{<-}{#1}{#2}}
 \rput[B](0.5,0){\pstree[treemode=R,levelsep=0.5,treesep=0.4]{\Tc*{0.6mm}}{
     \pstree{\Tc*{0.6mm}}{\Tc*{0.6mm}
     \Tc*{0.6mm}}  \pstree{\Tc*{0.6mm}}{\Tc*{0.6mm} \Tc*{0.6mm}}}}
 \end{pspicture}
 \caption{On the left: a tree for which $(\sigma,<,T')$ is not an intermediate triple.
 Only vertices belonging to the tree $U$ were shown.
 On the right: the same tree with the opposite orientation of the edge $BD$. In this
 case the pair of edges $BA$, $DB$ does not form a bay. Order of vertices is given by
 $R\leq A<B<D$.}
 \label{fig:inversetransform2bis}
\end{figure}
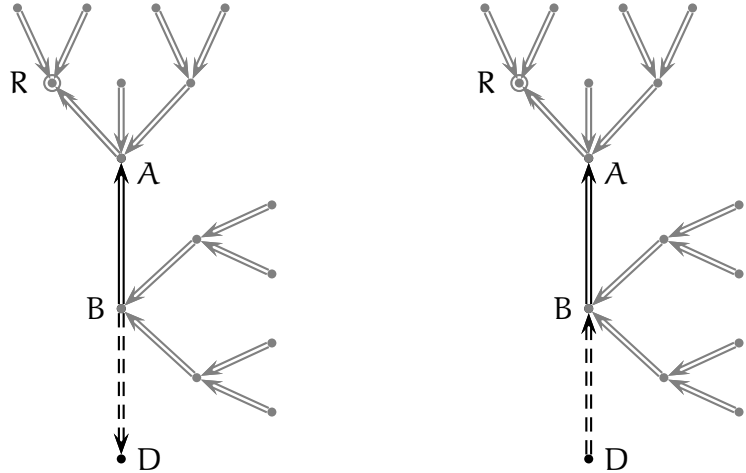

Let $U$ denote the subtree of $T_{\sigma}$ which consists of the
vertices $\{D\}\cup\{x:R \preceq x \text{ and } x\leq S\}$. Let us
change for a moment the orientation of the edge $DB$, as shown on
the right-hand side of Figure \ref{fig:inversetransform2} and the
right-hand side of Figure \ref{fig:inversetransform2bis}. It is easy
to see that after this change the tree $U$ has the form as
considered in Section \ref{subsec:leafsandbays}, i.e.\ the arrows on
all edges are pointing towards the root and it has at least two
vertices.

Let us consider the case when a corner formed in the vertex $B$ by
some edge on the left and the edge $DB$ on the right is a bay (on
the right-hand side of Figure \ref{fig:inversetransform2} this
corner corresponds to the pair of edges $uB$, $DB$). In this case
let $C$ denote the leaf corresponding to this bay, as described in
Section \ref{subsec:leafsandbays}. We unglue the edges $BA$ and $BD$
and we reglue them in the other way, as we described in Section
\ref{subsec:forward}, and thus we obtain a tree $T_{\sigma'}$
corresponding to some pairing $\sigma'$. Figure
\ref{fig:inversetransform0} describes the way how the vertices of
the original tree $T_{\sigma}$ are identified with the vertices of
$T_{\sigma'}$. We shall prove in the following that $(\sigma',<,S)$
is an intermediate triple.

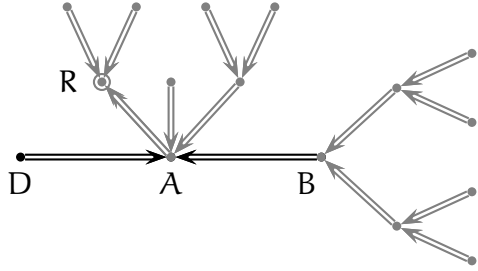
\begin{figure}
 \psset{unit=2cm}
 \begin{pspicture}[](-1,-0.8)(1,1.2)
 \cnode*(0,0){0.6mm}{A}
 \cnode*(1,0){0.6mm}{C}
 \cnode*(-1,0){0.6mm}{B}
 \ncline[arrowsize=2mm,doubleline=true]{<-}{A}{B}
 \ncline[arrowsize=2mm,doubleline=true]{<-}{A}{C}
 \rput[t](0,-0.1){$A$}
 \rput[t](0.9,-0.1){$B$}
 \rput[t](-1,-0.1){$D$}

\psset{linecolor=gray}

\renewcommand{\psedge}[2]{\ncline[doubleline=true,arrowsize=2mm]{->}{#1}{#2}}
 \rput[B](0,0){\pstree[treemode=U,levelsep=0.5,treesep=0.4]{\Tc*{0.6mm}}{
     \pstree{ \Tc{1.2mm}~[tnpos=l,tnsep=2mm]{$R$} \Tc*{0.6mm}  }{
     \renewcommand{\psedge}[2]{\ncline[doubleline=true,arrowsize=2mm]{<-}{#1}{#2}}
   \Tc*{0.6mm}  \Tc*{0.6mm}}
   \renewcommand{\psedge}[2]{\ncline[doubleline=true,arrowsize=2mm]{<-}{#1}{#2}}
   \Tc*{0.6mm}\pstree{\Tc*{0.6mm}}{\Tc*{0.6mm} \Tc*{0.6mm}}}}

  \renewcommand{\psedge}[2]{\ncline[doubleline=true,arrowsize=2mm]{<-}{#1}{#2}}
 \rput[B](1.5,0){\pstree[treemode=R,levelsep=0.5,treesep=0.4]{\Tc*{0.6mm}}{
     \pstree{\Tc*{0.6mm}}{\Tc*{0.6mm}
     \Tc*{0.6mm}} \pstree{\Tc*{0.6mm}}{\Tc*{0.6mm} \Tc*{0.6mm}}}}

 \end{pspicture}
 \caption{The tree from the left-hand side of Figure \ref{fig:inversetransform2bis} after ungluing
 and regluing differently the edges $BA$ and $BD$. Please note that the labels carried by
 the vertices $A,B,D$ have changed.}
 \label{fig:inversetransform0bis}
\end{figure}

It remains to consider the case when there is no bay in the vertex
$B$ formed by some edge on the left and the edge $DB$ on the right,
cf Figure \ref{fig:inversetransform2bis}. We again unglue and reglue
differently the edges $BA$ and $BD$; the resulting tree is denoted
by $T_{\sigma'}$. The identification of the vertices of $T_{\sigma}$
and $T_{\sigma'}$ is presented on Figures
\ref{fig:inversetransform2bis} and \ref{fig:inversetransform0bis};
please note that only vertices $A,B,D$ are nontrivially identified.

\begin{figure}
\psset{unit=2cm}
\begin{pspicture}[](-1,-1.1)(1,2)
 \cnode*(0,1){0.6mm}{A}
 \cnode*(0,0){0.6mm}{B}
 \cnode*(0,-1){0.6mm}{S}

 \ncline[arrowsize=2mm,doubleline=true]{<-}{A}{B}
 \ncline[arrowsize=2mm,doubleline=true,linestyle=dashed]{<-}{S}{B}
 \rput[l](0.1,0.9){$A$}
 \rput[l](0.2,0){$B$}
 \rput[l](0.1,-1){$D$}

\psset{linecolor=gray}
\renewcommand{\psedge}[2]{\ncline[doubleline=true,arrowsize=2mm]{->}{#1}{#2}}
 \rput[B](0,1){\pstree[treemode=U,levelsep=0.5,treesep=0.4]{\Tc*{0.6mm}}{
     \pstree{ \Tc{1.2mm}~[tnpos=l,tnsep=2mm]{$R$} \Tc*{0.6mm}  }{
     \renewcommand{\psedge}[2]{\ncline[doubleline=true,arrowsize=2mm]{<-}{#1}{#2}}
   \Tc*{0.6mm}  \Tc*{0.6mm}}
   \renewcommand{\psedge}[2]{\ncline[doubleline=true,arrowsize=2mm]{<-}{#1}{#2}}
   \Tc*{0.6mm}\pstree{\Tc*{0.6mm}}{\Tc*{0.6mm} \Tc*{0.6mm}}}}

 \renewcommand{\psedge}[2]{\ncline[doubleline=true,arrowsize=2mm]{<-}{#1}{#2}}
 \rput[B](-0.55,0){\pstree[treemode=L,levelsep=0.5,treesep=0.4]{\Tc*{0.6mm}}{
     \pstree{\Tc*{0.6mm}}{\Tc*{0.6mm} \Tc*{0.6mm}} \Tc*{0.6mm}\pstree{\Tc*{0.6mm}}{\Tc*{0.6mm}
     \Tc*{0.6mm}~[tnpos=l,tnsep=0mm]{$C$}    }}}

 \rput[B](0.5,0){\pstree[treemode=R,levelsep=0.5,treesep=0.4]{\Tc*{0.6mm}}{
     \pstree{\Tc*{0.6mm}}{\Tc*{0.6mm}
     \Tc*{0.6mm}}  \pstree{\Tc*{0.6mm}}{\Tc*{0.6mm} \Tc*{0.6mm}}}}
 \end{pspicture}
\hspace{2cm}
 \begin{pspicture}[](-1,-1.1)(1,2)
 \cnode*(0,1){0.6mm}{A}
 \cnode*(0,0){0.6mm}{B}
 \cnode*(0,-1){0.6mm}{S}

 \ncline[arrowsize=2mm,doubleline=true]{<-}{A}{B}
 \ncline[arrowsize=2mm,doubleline=true,linestyle=dashed]{->}{S}{B}
 \rput[l](0.1,0.9){$A$}
 \rput[l](0.2,0){$B$}
 \rput[l](0.1,-1){$D$}

\psset{linecolor=gray}
\renewcommand{\psedge}[2]{\ncline[doubleline=true,arrowsize=2mm]{->}{#1}{#2}}
 \rput[B](0,1){\pstree[treemode=U,levelsep=0.5,treesep=0.4]{\Tc*{0.6mm}}{
     \pstree{ \Tc{1.2mm}~[tnpos=l,tnsep=2mm]{$R$} \Tc*{0.6mm}  }{
     \renewcommand{\psedge}[2]{\ncline[doubleline=true,arrowsize=2mm]{<-}{#1}{#2}}
   \Tc*{0.6mm}  \Tc*{0.6mm}}
   \renewcommand{\psedge}[2]{\ncline[doubleline=true,arrowsize=2mm]{<-}{#1}{#2}}
   \Tc*{0.6mm}\pstree{\Tc*{0.6mm}}{\Tc*{0.6mm} \Tc*{0.6mm}}}}

 \renewcommand{\psedge}[2]{\ncline[doubleline=true,arrowsize=2mm]{<-}{#1}{#2}}
 \rput[B](-0.55,0){\pstree[treemode=L,levelsep=0.5,treesep=0.4]{\Tc*{0.6mm}}{
     \pstree{\Tc*{0.6mm}}{\Tc*{0.6mm} \Tc*{0.6mm}}
     \Tc*{0.6mm}\pstree{\Tc*{0.6mm}~[tnpos=b,tnsep=1mm]{$u$}}{\Tc*{0.6mm}
     \Tc*{0.6mm}~[tnpos=l,tnsep=0mm]{$C$}    }}}

 \rput[B](0.5,0){\pstree[treemode=R,levelsep=0.5,treesep=0.4]{\Tc*{0.6mm}}{
     \pstree{\Tc*{0.6mm}}{\Tc*{0.6mm}
     \Tc*{0.6mm}}  \pstree{\Tc*{0.6mm}}{\Tc*{0.6mm} \Tc*{0.6mm}}}}
 \end{pspicture}

 \caption{On the left: a tree for which $(\sigma,<,S')$ is not an intermediate
triple.
 Only vertices belonging to the tree $U$ were shown.
 On the right: the same tree with the opposite orientation of the edge $DB$. In this
 case the pair of edges $uB$, $DB$ forms a bay. The order of the vertices is given by
 $R\leq A<B<C\leq S<D$.}
 \label{fig:inversetransform2}
\end{figure}
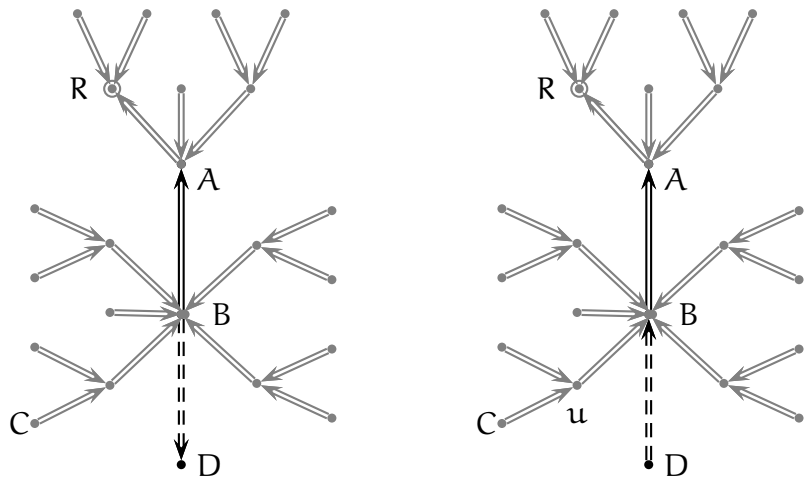

 \begin{figure}
 \psset{unit=2cm}
 \begin{pspicture}[](-1,-0.8)(1,1.2)
 \cnode*(0,0){0.6mm}{A}
 \cnode*(1,0){0.6mm}{C}
 \cnode*(-1,0){0.6mm}{B}
 \ncline[arrowsize=2mm,doubleline=true]{<-}{A}{B}
 \ncline[arrowsize=2mm,doubleline=true]{<-}{A}{C}
 \rput[t](0,-0.1){$A$}
 \rput[t](0.9,-0.1){$C$}
 \rput[t](-0.9,-0.1){$B$}

\psset{linecolor=gray}

\renewcommand{\psedge}[2]{\ncline[doubleline=true,arrowsize=2mm]{->}{#1}{#2}}
 \rput[B](0,0){\pstree[treemode=U,levelsep=0.5,treesep=0.4]{\Tc*{0.6mm}}{
     \pstree{ \Tc{1.2mm}~[tnpos=l,tnsep=2mm]{$R$} \Tc*{0.6mm}  }{
     \renewcommand{\psedge}[2]{\ncline[doubleline=true,arrowsize=2mm]{<-}{#1}{#2}}
   \Tc*{0.6mm}  \Tc*{0.6mm}}
   \renewcommand{\psedge}[2]{\ncline[doubleline=true,arrowsize=2mm]{<-}{#1}{#2}}
   \Tc*{0.6mm}\pstree{\Tc*{0.6mm}}{\Tc*{0.6mm} \Tc*{0.6mm}}}}

 \renewcommand{\psedge}[2]{\ncline[doubleline=true,arrowsize=2mm]{<-}{#1}{#2}}
 \rput[B](-1.6,0){\pstree[treemode=L,levelsep=0.5,treesep=0.4]{\Tc*{0.6mm}}{
     \pstree{\Tc*{0.6mm}}{\Tc*{0.6mm} \Tc*{0.6mm}} \Tc*{0.6mm}\pstree{\Tc*{0.6mm}}{\Tc*{0.6mm}
     \Tc*{0.6mm}~[tnpos=l,tnsep=0mm]{$D$}    }}}

 \rput[B](1.5,0){\pstree[treemode=R,levelsep=0.5,treesep=0.4]{\Tc*{0.6mm}}{
     \pstree{\Tc*{0.6mm}}{\Tc*{0.6mm}
     \Tc*{0.6mm}} \pstree{\Tc*{0.6mm}}{\Tc*{0.6mm} \Tc*{0.6mm}}}}

 \end{pspicture}
 \caption{The tree from the left-hand side of Figure \ref{fig:inversetransform2} after ungluing
 and regluing differently edges $BA$ and $BD$. Please note that the labels carried by the vertices
 $A,B,C,D$ have changed.}
 \label{fig:inversetransform0}
\end{figure}
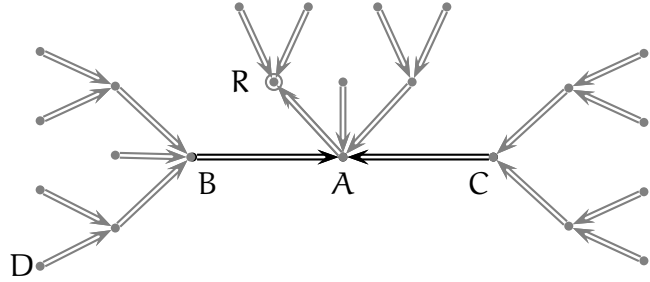

\begin{lemma}
The triple $(\sigma',<,S)$ given by the above construction is an
intermediate triple.

The above procedure will stop after a finite number of steps.
\end{lemma}
\begin{proof}
In the following we consider only the case presented on Figure
\ref{fig:inversetransform2} since the other one is analogous.

Conditions (\ref{interm1}) and (\ref{interm2}) are very easy to
verify. Condition (\ref{interm3}) holds true for tree $T_{\sigma}$
hence there are only two reasons why it could fail for the tree
$T_{\sigma'}$. Firstly, there might be some pair of adjacent
vertices $v,w$ on the tree $T_{\sigma'}$ for which $v\prec w$ and
$v>w$ which is new, i.e.\ which is not present in the tree
$T_{\sigma}$. There are three possible cases: $v=D$, $C<w<D$
(impossible, since \eqref{equ:nicmiedzyTiD} implies $C<w\leq S<D$
and $w\in U$ which contradicts that $C$ is a leaf of the tree $U$);
or $A<v<D$, $w=A$ (causes no problems); or $v=C$, $B<w<C$ (causes no
problems).

The second reason why condition (\ref{interm3}) could fail is that
for some pair of adjacent vertices $v,w\notin\{A,B,C,D\}$ such that
$v \prec w$ and $v>w$ the set
\begin{equation}
\label{eq:nibyniepusty} \{x: R\preceq x\text{ and } S<x<v\}
\end{equation}
might be non-empty. Since tree $T_{\sigma}$ fulfils condition
(\ref{interm3}) therefore the only element which could possibly
belong to the set \eqref{eq:nibyniepusty} is $D$. This, however, is
not the case since $D$ was chosen to be the maximal element in the
set of the possible values of $v$.

Note that in each step of our operation the cardinality of the set
$\{x\in T_{\sigma}: R\preceq x \text{ and } S\leq x\}$ increases.
This shows that our procedure will eventually stop.
\end{proof}

The proof of the following lemma is straightforward and we leave it
to the Reader.
\begin{lemma}
The operation described in Section \ref{subsec:forward} and the
operation described in Section \ref{subsec:backward} are inverses of
each other.
\end{lemma}

Thus, the proof of Theorem \ref{theo:malabijekcja} is finished.


We will find the following result useful in Section
\ref{sec:themainresult}.
\begin{lemma}
\label{lem:inclusion} Let $(\sigma,<,S)$, $(\sigma',<,S')$ be
intermediate triples such that $(\sigma',<,S')$ is obtained from
$(\sigma,<,S)$ by a number of forward transformations from Section
\ref{subsec:forward}. We denote $ U=\{v\in T_{\sigma}: R\preceq v
\text{ and } v\leq S \} $ and $U'= \{v\in T_{\sigma'}: R\preceq v
\text{ and } v\leq S' \}$.

Then $U\subseteq U'$. Furthermore, every element of the difference
$U'\setminus U$ is bigger (both with respect to the order $<$ and
$\lhd$) than every element of the set $U$.

Secondly,
\begin{equation}
\label{eq:zawieranie}
\{v\in T_{\sigma}: R\preceq v \}\supseteq \{v\in T_{\sigma'}:
R\preceq v \}.
\end{equation}

\end{lemma}
Inductive proof is straightforward. Please note that, contrary to
the order $<$, the order $\lhd$ is different on the trees
$T_{\sigma}$ and $T_{\sigma'}$ and the lemma holds true for both
choices of $\lhd$.

\section{Proof of the correctness of the main bijection}

\label{sec:themainresult}

We will prove Theorem \ref{theo:mainresult}, namely that \MainBijection indeed
provides the desired bijection.


%
%
%
%
%
%
%
%
%

\subsection{Intermediate points}
\label{subsec:intermediate}

Our strategy is to describe possible intermediate outcomes of the algorithm
\MainBijection. 

\begin{definition}
We call $\big(T_{\sigma},  (B_{i+1},\dots,B_{m}) \big)$ an
intermediate point for $i$ ($0 \leq i\leq m$) if

\begin{enumerate}

\item $\sigma$ is a pairing compatible with $\epsilon_i$, as defined in
\eqref{eq:epsilon};

\item $T_{\sigma}$ is a quotient tree with some of the vertices labeled with
different elements of $\{1,\dots,L\}$;

\item let $V$ be the set of unlabeled vertices of $T_{\sigma}$; if $i=m$ then
$V=\emptyset$, if $i<m$ then $V$ is a tree such that $R\in V$, $V\subseteq\{x\in
T_{\sigma}: x \succeq R\}$, furthermore for all pairs such that $x\prec y$ and
$y\in V$ we also have $x\in V$; \label{intermediate3}

\item for all pairs of labeled vertices such that $x\prec y$ their labels
fulfill $x<y$; \label{intermediate4}

\item sets $B_{i+1},\dots,B_{m}$ and the set of labels are disjoint and their
union is equal to $\{1,\dots,L\}$;

\item $|B_n|+\cdots+|B_m|\geq l_n+\cdots+l_m+1$ for all $i+1\leq n\leq m$.



%
%

\end{enumerate}
\end{definition}

\subsection{Startpoints and endpoints}
\begin{lemma}
There is a bijection between the elements of the set {\MainPary} (thus input
data of algorithm \MainBijection) and
the intermediate points corresponding to $i=m$.

There is a bijection between the elements of the set {\MainZbiory}
and the intermediate points corresponding to $i=0$ (for which algorithm
\MainBijection terminates).
\end{lemma}

%

\subsection{The forward transformation}
\label{subsec:forward2}

\begin{lemma}
\label{lem:MainBijection-forward-ok}
After each iteration of the loop in \MainBijection tuple
$\big(\tilde{T},(B_i,\dots,B_m) \big)$ forms an intermediate point, where
$\tilde{T}$ denotes the tree $T$ with all artificial labels removed.
Futhermore, on the set of the vertices with artificial labels the order $\lhd$
coincides with the order of the labels $<$.

%
%
%
%
%
%
%
%
%
%
\end{lemma}
\begin{proof}
We are going to use backward induction with respect to $i$.

Firstly, observe that 
in line \ref{linia:malabijekcja} of \MainBijection the computation of
\SmallBijection{T} is performed on a tree $T$ which is as prescribed in point
(A) of Theorem \ref{theo:malabijekcja} therefore
afterwards $T$ is as presecribed in point (B) of Theorem
\ref{theo:malabijekcja}.

Furthermore, Lemma \ref{lem:inclusion} shows that in line \ref{lin:remove} of
\MainBijection all artificial labels will be removed. Therefore all artificial
labels in $T$ (equivalently, unlabeled vertices $\tilde{T}$) after the
iteration of the loop must have been created in line \ref{linia:koniec}.

This shows points (\ref{intermediate3}) and (\ref{intermediate4}) in the
definition of an intermediate point.
\end{proof}


%

\subsection{The backward transform}
\label{subsec:backward2} We shall prescribe now the inverse of the
transform prescribed in Section \ref{subsec:forward2} (i.e.\ a single
interation of the loop in \MainBijection). Since we
simply have to reverse all steps of the forward transformation, our
description will be quite brief and we shall concentrate only on the
most critical points.

As we pointed out in the proof of Lemma \ref{lem:MainBijection-forward-ok} all
artificial labels in $T$ (equivalently, unlabeled vertices $\tilde{T}$) after
the iteration of the loop must have been created in line \ref{linia:koniec}; in
other words, tree $U$ consists of vertices carrying artificial labels.

Therefore, in order to undo line \ref{linia:koniec} of \MainBijection we simply
remove all artificial labels and in order to undo line
\ref{lin:catalanglue} we unglue all edges with both unlabeled ends.

In order to undo line \ref{linia:flip} we change the orientation of all edges
and reverse the order of $<$. In order to undo line \ref{linia:remove-edges} to
both sides of the root we attach $l_i$ new edges with appropriate orientations.

In order to undo line \ref{linia:drugi-krok} we glue the unpaired edges
according to the Catalan pairing.

Let $U$ denote the set of unlabeled vertices.
In order to undo line \ref{lin:remove}
we create $|U|-|B_i|$ artificial labels (which are integer numbers
different from $1,\dots,L$) which are smaller than any element of the set
$\{1,\dots,L\}$ and we label the elements of $U$ with these artificial
labels and the labels from the set $B_i$ in such a way that the order
$\lhd$ of the vertices of $U$ coincides with the order $<$ on the labels.

In this way tree $T$ is as prescribed in point (B) of Theorem
\ref{theo:malabijekcja} hence undoing line \ref{linia:malabijekcja} of
\MainBijection is possible.

%
%
%

\begin{lemma}
If the above procedure was started for an intermediate point
$\big(\tilde{T}^{\text{initial}}, (B_{i},\dots,B_m) \big)$ for $i'=i-1$ then
the resulting tuple $\big(\tilde{T}^{\text{final}}, (B_{i+1},\dots,B_m) \big)$
is an intermediate point for $i$, where $\tilde{T}$ is the tree $T$ given by the
above procedure with all artificial labels removed.
\end{lemma}
\begin{proof}
Only condition (\ref{intermediate3}) is less trivial and needs
to be proved. The inclusion \eqref{eq:zawieranie} shows that after
undoing of line \ref{linia:malabijekcja} every unlabeled vertex $x$ fulfills
$x\succeq R$. Furthermore, in tree $T$ the order of the labels
is compatible with the orientations of the edges; therefore if $x\prec y$ and
$y$ is unlabeled then $x<y$ and $x$ is also unlabeled.
\end{proof}

We leave the following lemma as a simple exercise to the Reader.

\begin{lemma}
Operations described in Sections \ref{subsec:forward2} and
\ref{subsec:backward2} are inverses to each other.
\end{lemma}
Thus, the proof of Theorem \ref{theo:mainresult} is finished.

%
%
%
%

%
%
%
%
%
%



\section{Acknowledgments}

I thank Kenneth Dykema for introducing me into the subject.
Research supported by the MNiSW research grant P03A 013 30, by the EU Research
Training Network ``QP-Applications", contract HPRN-CT-2002-00279 and  by the EC
Marie Curie Host Fellowship for the Transfer of Knowledge ``Harmonic Analysis,
Nonlinear Analysis and Probability", contract MTKD-CT-2004-013389.
The author is a holder of a scholarship of European Post-Doctoral
Institute for Mathematical Sciences.

%
%

\end{document}